\documentclass[twoside,11pt]{amsart}

\newcommand{\ra}{\rightarrow}

\newcommand{\BR}{\mathbb{R}}
\newcommand{\CP}{\mathbb{CP}}

\newcommand{\BZ}{\mathbb{Z}}

\newcommand{\Kahler}{K\"{a}hler }

\newcommand{\M}{\mathfrak{M}}

\newcommand{\TryPackage}[3]{\IfFileExists{#1.sty}{\usepackage{#1}#2}{#3}
}
\TryPackage{mathrsfs}{\renewcommand{\mathcal}{\mathscr}}{%
        \TryPackage{eucal}{}{}}

\newcommand{\ga}{\gamma}

\newcommand{\ZZ}{{\mathbb Z}}
\newcommand{\RR}{{\mathbb R}}
\newcommand{\CC}{{\mathbb C}}

\newcommand{\cM}{{\mathcal M}}

\usepackage{latexsym}
\usepackage{amsmath}
\usepackage{amsthm}
\usepackage{amscd}
\usepackage{xypic} 
\usepackage{amsfonts}
\usepackage{amssymb}
\usepackage{graphicx}
\usepackage{psfrag}

\newtheorem{df}{Definition}
\newtheorem{thm}[df]{Theorem}
\newtheorem{cor}[df]{Corollary}
\newtheorem{lem}[df]{Lemma}
\newtheorem{prop}[df]{Proposition}
\newtheorem{ex}[df]{Example}

\newtheorem{conj}[df]{Conjecture}
\newtheorem{quest}{Question}

\begin{document}

\title{On symplectic 4-manifolds with prescribed fundamental group}

\author{Scott Baldridge}
\author{Paul Kirk}

\thanks{The first   author  gratefully acknowledges support from the
NSF  grant DMS-0507857. The second  author  gratefully acknowledges
support from the NSF  grant DMS-0202148.}

\address{Department of Mathematics, Louisiana State University \newline
\hspace*{.375in} Baton Rouge, LA 70817}
\email{\rm{sbaldrid@math.lsu.edu}}

\address{Mathematics Department, Indiana University \newline
\hspace*{.375in} Bloomington, IN 47405}
\email{\rm{pkirk@indiana.edu}}

\maketitle

\begin{abstract} In this article we study the problem of minimizing
$a\chi+b\sigma$
on the class of all symplectic 4--manifolds with prescribed
fundamental group $G$ ($\chi$ is the Euler characteristic, $\sigma$
is the signature, and $a,b\in \BR$), focusing on the important cases $\chi$, $\chi+\sigma$ and $2\chi+3\sigma$.   In certain
situations we can derive lower bounds for these functions and
describe symplectic 4-manifolds which are minimizers.  We
derive an upper bound for the minimum of  $\chi$ and $\chi+\sigma$
in terms of the presentation of $G$.

\end{abstract}


\section{Introduction}

Pick a finitely presented group $G$ and  let $\M(G)$ denote the  class
of  closed symplectic 4-manifolds $M$ which have $\pi_1(M)$ isomorphic
to $G$. The existence of a
symplectic $M$  with given fundamental group $G$ was demonstrated
by Gompf \cite{symp:gompf:construct_symp_man}.

In this article we study the  problem of finding  minimizers  in
$\M(G)$ where minimizing is taken  with regard to the Euler
characteristic
$\chi$, following the approach introduced by Hausmann and Weinberger in
\cite{group:euler_character_of_fund_group:hausmann} for smooth
4-manifolds.    There are two aspects to this problem. Finding lower
bounds to $\chi(M)$ for $M\in\M(G)$ addresses the question ``How
large must a symplectic manifold with fundamental group $G$ be?''
The other aspect of the problem is finding efficient and explicit
constructions of symplectic manifolds with a given fundamental
group.

Our main general result concerning upper bounds is Theorem \ref{thm1},
which states:

\medskip

\noindent{\bf Theorem  \ref{thm1}.} {\em Let  $G$ have a   presentation
with $g$ generators
$x_1,\cdots, x_g$ and $r$ relations $w_1,\cdots,w_r$. Then
   there exists a closed symplectic 4-manifold $M$ with $\pi_1M\cong G$,
Euler characteristic
$\chi(M)=12(g+r+1)$, and signature $\sigma(M)=-8(g+r+1)$.}

\medskip

We also provide a number of examples of small closed symplectic
manifolds with certain fundamental groups.  A successful example is the
following theorem,  which generalizes to the symplectic  setting the
results of
\cite{group:h-w_4-man_invar_abel_grp:kirk}. (See Corollary
\ref{evenZn}, for a more complete statement.)

\medskip

\noindent{\bf Theorem.} {\em Let $F_g$ denote the  closed oriented
surface of genus $g$, and let
$S_g=${\rm Sym}$^2(F_g)$, so  that $S_g$ is a closed symplectic
manifold with  fundamental group $\ZZ^{2g}$. If $g\equiv 0,1,$ or $3
\pmod{4}$, then any other  closed symplectic 4-manifold
$N$ with $\pi_1(N)\cong \ZZ^{2g}$ satisfies $\chi(N)\ge \chi(S_g)$.}

\medskip

The general theme of this article is to investigate the simplest
symplectic 4-manifolds one  fundamental group at  a  time, finding
constructions, obstructions, and examples of minimizers of   $a\chi+b\sigma$.

The problem of   minimizing $\chi$ and $\chi+\sigma$  of
4-manifolds with a prescribed  fundamental group arises in many contexts and has  been studied in  explicitly in a number of interesting articles.
Hausmann and Weinberger in \cite{group:euler_character_of_fund_group:hausmann}
used $q(G)=\min_{\pi_1(M^4)\cong G}\chi(M)$ to establish the existence of a perfect group  which can be the fundamental group of
a homology sphere in dimensions greater than 4 but which is not the fundamental group of a homology 4-sphere, and to construct groups  which are knot groups in  dimensions greater than 4 but which are not the fundamental group of a knotted 2-sphere in $S^4$.

Kotschick in \cite{kot2}  inserted the signature into the topic by defining the invariant $p(G)=\min_{\pi_1(M^4)\cong G}\chi(M)-|\sigma(M)|$  and in \cite{kot1}  he carries out a systematic study of
$p(G)$ and $q(G)$, including computations and estimates for  $q(G)$ and $p(G)$ for various $G$.
Moreover, Kotschick  discusses the problem of
defining variants of $p$ and $q$ by restricting to 4-manifolds with fundamental group  $G$ which admit various geometric structures, e.g.
spin structures, almost complex structures, positive scalar curvature, and, symplectic structures, the topic of the present article.  He also investigates the question of what the possible  values of $p(G)$ and  $q(G)$ are for a given group $G$, a question that we generalize and recast in Section  3.

Other  related work   includes the articles of Eckmann  \cite{eck}  and L\"uck \cite{luck} who derive  bounds on $p(G)$ and $q(G)$ for various $G$ using $\ell^2$-cohomology and the $\ell^2$-signature theorem, as well as the articles \cite{bohr}, \cite{KJ}, \cite{group:h-w_4-man_invar_abel_grp:kirk}. The
general problem of calculating $q(G)$ appears as Problem 4.59 of Kirby's problem list \cite{kirp}.

\medskip

The article  is organized as follows.  In Section 2 we establish some
simple bounds and describe Gompf's construction for producing a
symplectic 4--manifold with a given fundamental group.  The function
$f=a\chi+b\sigma$ for $a,b\in \BR$ is studied in Section 3 and some
reasons are given for restricting to the cases $\chi$ and
$\chi+\sigma$.  In Section 4 we describe new constructions that give
upper bounds for $\min \chi$ and $\min \chi+\sigma$ based upon the
group presentation of $G$.  In Section 5 we focus on examples for
specific classes of groups, namely free groups, cyclic groups, and
free abelian groups and describe minimizers of $\chi$ for many free
abelian groups.  In the last section, we  speculate about when or
whether there are conditions for which the minimizers of $\chi$ or
$\chi+\sigma$ are unique.

\bigskip

The authors would like to thank R. Gompf  and D. Kotschick for making helpful
and insightful comments which improved this article.

\section{Some bounds}
The fundamental  numerical  invariants of a  4-manifold are its
Euler  characteristic $\chi$ and its signature $\sigma$.  We will
focus on  the problem of minimizing $\chi$ and sometimes
$\chi+\sigma$ over the collection  of symplectic  manifolds with
fundamental group $G$.   Section \ref{sec:restrict} gives partial
justification for our restricting to these cases.
   We  remind  the reader of some coarse bounds on the
Euler characteristic of {\em smooth} closed orientable 4-manifolds
introduced in \cite{group:euler_character_of_fund_group:hausmann}.
First recall   that if $G$ is finitely generated and $M$ is a
connected oriented 4-dimensional Poincar\'e complex then the second
Betti number of $M$, $b_2(M)$, is at least as large as the second
Betti number of $K(G,1)$ (with any field coefficients). Since
$b_1(M)=b_3(M)=b_1(G)$,  this implies:
\begin{equation}\label{hopf}
2-2b_1(G)+b_2(G)\leq \chi(M).
\end{equation}

By taking the double of the 2-handlebody defined by a presentation
of $G$, one obtains a smooth manifold $M$ with $\pi_1(M)=G$ and
$\chi(M)=2-2d$ where $d$ denotes the deficiency   of the
presentation (i.e.~the number of generators minus the number of
relations). Thus one has the bound {\em for smooth manifolds}, where
def$(G)$ denotes the minimum of the deficiency over all presentations:
\begin{equation*}\label{deficiency}
\min_{\pi_1M\cong G} \chi(M) \leq 2-2\text{def}(G).
\end{equation*}
This construction does not give a symplectic manifold in  general.
Thus this upper bound need not hold when one minimizes over
symplectic manifolds  with fundamental  group $G$. To obtain a
similarly general upper bound  requires an examination of the
construction of symplectic  manifolds  with  prescribed fundamental
group.

In the symplectic setting,  Gompf has given a  construction
\cite{symp:gompf:construct_symp_man}  by taking appropriate fiber sums
of $F\times T^2$
with many copies of the elliptic fibration $E(1)$. By examining
Gompf's argument one can formalize an upper bound.

Note that any finitely  presented  group is the quotient of
an oriented surface group, since (for example) the free group  on $g$
generators is a quotient  of the fundamental group  of a genus $g$
surface. Call a system of  immersed curves in general position
$\ga_i:S^1\to F, \ i=1,\cdots, r$ on  an orientable surface $F$ a
{\em geometric  surface presentation of $G$}  provided the
fundamental group of  the 2-complex obtained by attaching 2-cells
to $F$  along the $\ga_i$ is isomorphic to $G$.

Given a geometric surface presentation of  $G$,  the union of the
$\ga_i$
form  a  graph  $\Gamma$ (where one allows a graph to have some
isolated circle components).    Gompf's construction yields the
following general bound.

\begin{thm}[Gompf]\label{gompflemma}  Given   any  geometric surface
presentation  for
   $G$  with $r$ curves $\ga_1,\cdots, \ga_r$, if the  associated graph
$\Gamma$ has $n$ edges,
     there exists a  closed symplectic 4-manifold $M$ with $\pi_1(M)\cong
G$, $\chi(M)=12(r+2n+1)$ and $\sigma=-8(r+2n+1)$.
     Moreover, there exists a spin symplectic 4-manifold with
$\pi_1(M)\cong G$, $\chi(M)=24(r+2n+1)$ and $\sigma(M)=-16(r+2n+1)$.
\qed
   \end{thm}

Simple  experiments show that the number $n$ in Theorem
\ref{gompflemma} can be quite large for even simple group
presentations.  As an example we compute the Euler characteristic of a
manifold which has $G=\ZZ^4$. In this
situation, start with a genus 4 surface $F$ with
a standard collection of oriented circles
$$\alpha_1,\alpha_2,\alpha_3,\alpha_4,
\beta_1,\beta_2,\beta_3,\beta_4$$ in $F$ representing a symplectic
basis of $H_1(F)$.  The quotient
$\pi_1(F)/\langle\beta_1,\ldots,\beta_4\rangle$ is a free group
generated by the $\alpha_i$'s.  For $i=1,\ldots, 4$, let $\gamma_i
= \beta_i$. For $i=1,2,3$, set $\gamma_{i+4} = [\alpha_i,
\alpha_{i+1}]$ using the configuration of curves on the top of $F$
shown in Figure 1. 

\begin{figure}[h]\label{GET}
\psfrag{a}{$\alpha_1$} \psfrag{b}{$\alpha_2$}\psfrag{c}{$\alpha_3$}
\psfrag{d}{$\alpha_4$}
\begin{center}\includegraphics[scale=.7]{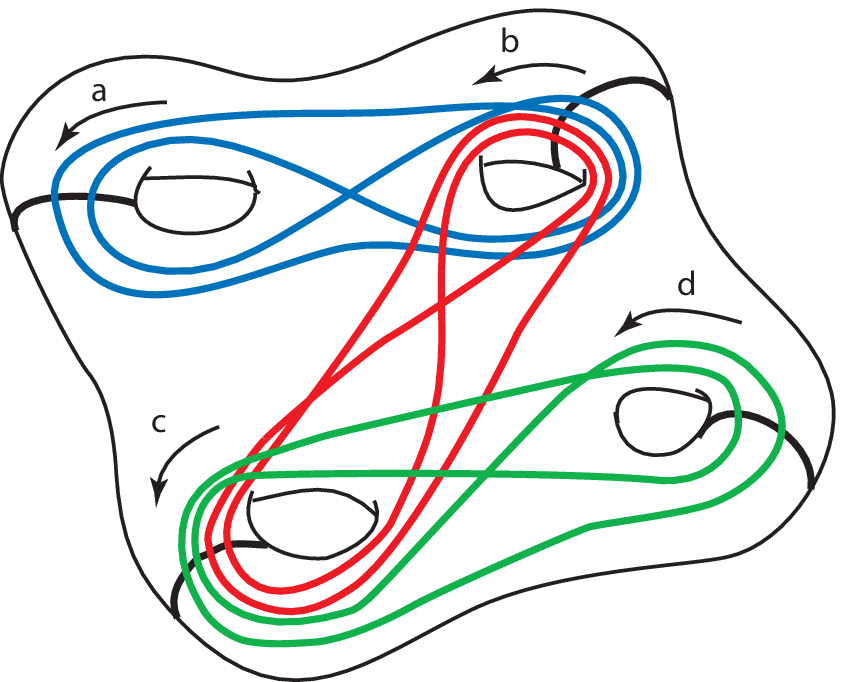}
\end{center} \caption{}
\end{figure}

\noindent Finally, set $\gamma_8 = [\alpha_2,\alpha_4]$,
$\gamma_9=[\alpha_1,\alpha_4]$, $\gamma_{10}=[\alpha_1,\alpha_3]$
using the same configuration as in Figure 1, but now on the
bottom of $F$, i.e. $\gamma_5,\gamma_6,\gamma_7$ are disjoint from
$\gamma_8,\gamma_9,\gamma_{10}$.  The union of the immersed curves
$\gamma_1,\ldots,
\gamma_{10}$ is an example of a geometric surface presentation of
$\ZZ^4$. After a careful count one finds 136 edges in the graph
described above.  Using the theorem above one computes:

\begin{ex}
The construction above produces a symplectic manifold $M$ with
$\pi_1(M)=\ZZ^4$ and $\chi(M)=3,396$.
\end{ex}

Gompf was not trying to minimize the Euler characteristic in his
construction; in fact, it is clear from his writings that he knows ways
to reduce this number significantly.  Still, our best estimate using
this construction as the starting point together with some tricks
(known to us) is $\chi(M)=516$.  This is a significant reduction, no
doubt, but the 4-torus $T^4$ has fundamental group $\ZZ^4$ and
$\chi=0$.  Thus
constructions like this one do not  give a particularly effective upper
bound for $\chi(M)$ for $M$ symplectic, $\pi_1(M)=G$. Moreover,
from the point of view of the present article the problem of
expressing a bound on the number $n$ of edges of $\Gamma$ in terms
of algebraic invariants of $G$ is unwieldy in general.

\bigskip

We end this section by recalling two facts that are
    useful in increasing the lower bound of Equation (\ref{hopf}) for
symplectic manifolds.  First, the symplectic form $\omega$ on a
symplectic 4-manifold $M$ has the property that $\omega\wedge
\omega$  is a volume form.  Thus $b^+(M)$, the dimension of the
largest  positive definite subspace of the intersection form (over
$\RR$) is always at least 1, and in particular, the second Betti
number $b_2(M)\ge 1$. For example,  this implies that if
$\pi_1(M)\cong\ZZ$, then $1\leq  \chi(M)$, improving Equation
(\ref{hopf})
by one when $G=\ZZ$.

Secondly, a symplectic  manifold  admits an almost complex
structure.   This has implications on its characteristic classes.
The consequence of most use to us is that $1-b_1(M)+ b^+(M)$ (the
index of the ASD complex) is even. For example, if $M$ is symplectic
and  $\pi_1(M)\cong\ZZ$, then $b^+(M)$  is even. Combined with the
observation of the previous paragraph, we conclude that $b^+(M)\ge
2$, and hence $2\leq\chi(M)$, improving Equation (\ref{hopf}) by two
when
$G=\ZZ$.

Putting these observations together one sees that if $M$ is symplectic with  fundamental
group $G$, then  $\chi(M)+\sigma(M)=2-2b_1(G) + 2b^+(M)$, and hence
\begin{equation}\label{chi+sig}
\chi(M)+\sigma(M)\ge \begin{cases} 4-2b_1(G)& \text {if } b_1(G) \text{
is even}, \\
6-2b_1(G)& \text {if } b_1(G) \text{ is odd.}
\end{cases}\end{equation}

\bigskip

\section{minimizing $a\chi +b \sigma$ and the special points $\chi$, $\chi +\sigma$,  and $2\chi
+3\sigma$}
\label{sec:restrict}

   In this section we investigate the values of $a$
and $b$ for which the function $a\chi+b\sigma$ has a lower bound on a suitable
class of 4--manifolds with a given fundamental group (smooth,
symplectic, etc.).  The answers to this question naturally  lead  to
breaking points  at  $a=b$,  and  $3a=2b$. These are related to important invariants of
symplectic 4--manifolds:
$\chi+\sigma$ is 4 times the holomorphic Euler characteristic, and
$2\chi +3\sigma$ is the square of the canonical class on a symplectic
manifold.  The approach described in this section can be viewed as a variant of the geography problem for 4-manifolds.

We first introduce a general  notion. Let $\M$ denote a class of closed oriented 4-manifolds. We will be most interested in the cases, $\M=\M(G)$, the class of symplectic 4-manifolds with fundamental group
$G$,   $\M=\M^\infty(G)$, the class of {\em smooth} manifolds with
fundamental group $G$, and
 $\M=\M^{min}(G)$, the subclass of $\M(G)$ consisting of minimal
symplectic 4-manifolds  with fundamental group  isomorphic  to $G$ (recall that a symplectic  4-manifold $M$ is called {\em  minimal} if  it is
not a blow up, i.e. $M\not \cong N\# \overline{\CC P}^2$ for $N$ symplectic). But the following result also applies in greater generality, e.g.    the class 4-dimensional Poincar\'e complexes with a given fundamental  group, or
  the class of smooth complex projective surfaces with a given fundamental group, or  the class of smooth 4-manifolds with even intersection form (for which the results of \cite{bohr} are relevant), or the class of almost complex 4-manifolds with given fundamental group (see \cite{kot3}), or even  the class of all topological oriented 4-manifolds (with no fundamental group restriction).

For $(a,b)\in \RR^2$, define $f_{\M}(a,b)\in \RR\cup \{\pm \infty\}$ to be the infimum
$$f_{\M}(a,b)=\inf_{M\in\M}\{ a\chi(M)+b\sigma(M)\},$$
with the understanding that $f_{\M}(a,b)=\infty$ if $\M$ is empty (e.g. if $\M$ is the class of \Kahler\ manifolds with fundamental group $\ZZ^3$). Define the {\em domain $D_\M$ of $\M$} to be the set
$$D_\M=\{(a,b) | \ f_\M(a,b)\ne-\infty\}.$$
Thus $D_\M$ is the set of $(a,b)$ so that $a\chi+b\sigma$ is bounded below on $\M$. Notice that $D_\M$ is a cone   since $f_\M(ra,rb)=rf_\M(a,b)$ when $r\ge 0$.  Furthermore, if $\M\subset \M'$ then $f_{\M'}(a,b)\leq f_{\M}(a,b)$, and hence $D_\M\supset D_{\M'} $.

Recall that a function $f$ on a convex set $S$ is {\em concave} if $f(tx+(1-t)y)\ge tf(x)+(1-t)f(y)$ for all $x,y\in S.$

\begin{thm} \label{convex}
The domain $D_\M$  is a   convex cone and
$f_\M$ is a continuous concave function on $D_\M$.
\end{thm}
\begin{proof} The proof is simple: each $M\in \M$ determines a half space $H_M\subset \RR^3$ by
$$H_M=\{ (a,b,c)\ | \ c\leq a\chi(M)+b\sigma(M)\}.$$
The intersection
$$I=\cap_{M\in \M}H_M$$
is a convex set whose projection to $\RR^2$ is  $D_\M$.  Thus $D_\M$ is convex.   Furthermore, if $(a,b)\in D_\M$, then
$f_\M(a,b)$ is the largest number $c$ so that $(a,b,c)\in I$; this is clearly continuous and concave.
\end{proof}

 Since $D_\M$ is a convex cone, it is either the entire plane (e.g. if $\M$ contains finitely many homotopy types) or else it is a cone with angle less than or equal to $\pi$.

 Interestingly, $D_\M$ need not be closed. For example, let $\M=\{M_k \}_{k=1}^ \infty$, where
$$M_k=2k^2\overline{\CC P}^2 \#(k^2-k) S^2\times S^2.$$
Then $\chi(M_k)=2+4k^2-2k$ and $\sigma(M_k)=-2k^2$. Thus
$a\chi+b\sigma$ has a lower bound if $2a>b$ or if $2a=b$ and $a\leq 0$. Otherwise, $a\chi+b\sigma$ is not bounded below. Thus $D_\M$ has cone angle $\pi$ which contains one of its boundary rays ($\{(-r,-2r)\ | r> 0\}$) but not the other ($\{(r,2r)\ | r> 0\}$.

\medskip

We focus now on the class $\M^\infty(G)$ of  smooth 4--manifolds with fundamental group $G$.  Blowing up (i.e. taking the connected sum with $\overline{\CP}^2$ )increases $\chi$ by 1  and decreases $\sigma$ by 1 without changing $G$. Thus $a\chi+b\sigma$ is not bounded below if $a-b<0$, and so $D_{\M^\infty(G)}$ is contained in the half-plane $\{ a\ge b\}$. Similarly, taking connected sums with  $\CP^2$ shows that $D_{\M^\infty(G)}$ is contained in the half-plane $\{ a\ge -b\}$.  Hence   $D_{\M^\infty(G)}$ lies in the cone $\{a\ge |b|\}$.

If $a\ge |b|$ then
$$
a\chi(M)+b\sigma(M)=2a(1-b_1(G)) + (a+b)b^+(M)+ (a-b)b^-(M)\ge 2a(1-b_1(G))
$$
and so $(a,b)\in D_{\M^\infty(G)}$. Thus we have proven the following.

\begin{prop} Fix a group $G$ and  $a\ne 0$. Then
$f_{\M^\infty(G)} $ has domain $$D_{\M^\infty(G)}=\{ (a,b)\ | \ a\ge |b| \}, $$
i.e. $D_{\M^\infty(G)}$ is the cone over the closed interval $\{1\}\times[-1,1]$.
\qed\label{prop:bound_below_M_smooth}
\end{prop}

Restricting to the class of symplectic manifolds $\M(G)$ everything
follows as above except for one point: taking connected sum of a
symplectic manifold with $\CP^2$ does not yield a symplectic manifold. In particular, one cannot conclude that  $a\chi+b\sigma$ has no lower bound
on $\M(G)$ for $a>-b$. Theorem 6.3 of \cite{symp:gompf:construct_symp_man} shows that there exists symplectic manifolds with
fundamental group $G$ and arbitrarily large signature. Thus $b\sigma$ does not have a
lower bound on $\M(G)$  when $b<0$.

These observations imply that
the domain $D_{\M(G)}$ is contained in the intersection of the    half-planes $b\leq a$  and $a\ge 0$,  and contains the ray $\{(r,r)\ | \ r\ge 0\}$ as one the two boundary edges of the cone $D_{\M(G)}$. The other  edge is a ray  $\{ (r \cos(\theta_G), r\sin(\theta_G)\ | \ r\ge 0 \}$ for  some angle $\theta_G$  in
$[-\tfrac{\pi}{2}, -\tfrac{\pi}{4}]$. We were unable to determine the ``critical'' angle $\theta_G$. This leads us to pose the question:

\begin{quest} Does the domain $D_{\M(G)}$ contain any pairs $(a,b)$ with $a>-b$?
Does   $\theta_G$ depend on the group $G$? \end{quest}

 For $G=\{e\}$, Stipsicz (\cite{stip}) has constructed simply connected symplectic 4-manifolds so  that $a\chi+b\sigma$ is not bounded below when $b<-\frac{10}{3}a$, so that
$\theta_{\{e\}}\ge\tan^{-1}(-\frac{10}{3})$.

Figure 2 explains the notation.
\begin{figure}\label{DM}
\includegraphics[scale=.7]{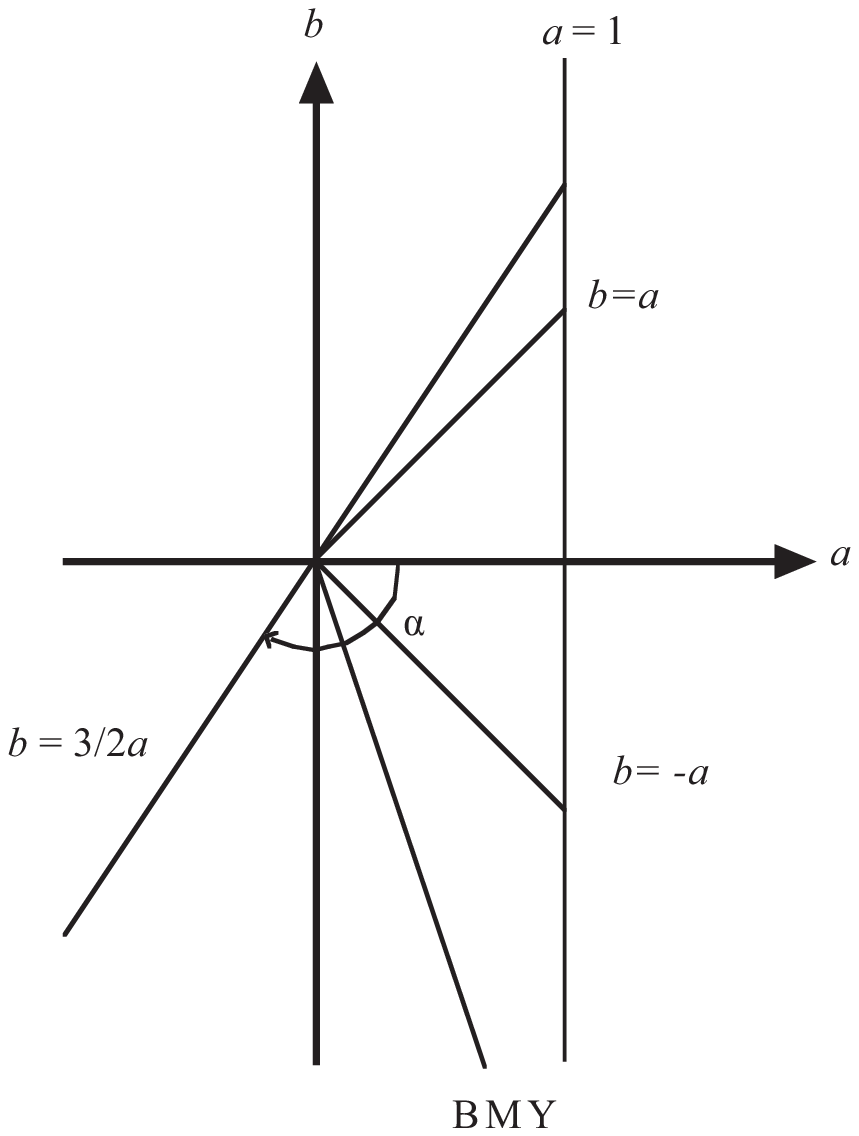}
\caption{}
\end{figure}
\medskip

We now look at the class $\M^{min}(G)$ of minimal  symplectic manifolds
with  fundamental group $G$.  This time blowing up is not allowed, since by definition minimal symplectic manifolds are not  blowups. Since $\M^{min}(G)\subset  \M^{\infty}(G)$
we know by Proposition \ref{prop:bound_below_M_smooth} that $D_{\M^{min}(G)}$  contains  the cone over the interval $\{1\}\times [-1,1]$.  The following proposition implies that $D_{\M^{min}(G)}$ is strictly larger than $D_{\M(G)}$.

\begin{prop} \label{prop:minimal_symp_man} Fix a group $G$.  Then
 $\chi+b\sigma$ has a lower bound on $\M^{min}(G)$ if $-1\leq
b\leq \frac32$ and does not if $b>\frac32$. In particular, $D_{\M^{min}(G)}$ (and hence $D_{\M(G)}$) is contained  in the half-plane $\{(a,b)\ | \ b\leq\frac32 a\}$, and
$D_{\M^{min}(G)}$ contains the cone over the interval $\{1\}\times [-1,\frac32]$.
\end{prop}

\begin{proof}
Let $K$ be the canonical class of $M\in \M^{min}(G)$.  A theorem of Liu
\cite{symp:app_gen_wall_cross} states that if $K^2<0$, then $M$ is
diffeomorphic to an irrationally ruled surface with fundamental group a
surface group.  Assume for a moment that $G$ is not a surface group.
In this case $K^2\geq 0$ or, equivalently $2\chi(M)+3\sigma(M)\ge 0$ for all
manifolds $M\in \M^{min}(G)$.  The convexity of the cone $D_{\M^{min}(G)}$ and the fact that $D_{\M^\infty(G)}\subset D_{\M^{min}(G)}$ implies that $D_{\M^{min}(G)}$
contains the cone $\{(a,b)\ | \ b\leq \frac32 a \text{ and } a\ge -b\}$.
The first part of the proposition
follows from this inequality for such groups.

The case when $G$ is a surface
group is similar. (Note that in this case there are only two manifolds in
$\M^{min}(G)$ up to diffeomorphism with $K^2 < 0$).

To prove that $\chi+b\sigma$ is unbounded when $b>\frac32$, let $M$ be
a spin symplectic manifold with $\pi_1M\cong G$ given by Gompf's
construction.  Then $ 2\chi(M)+3\sigma(M)=0$. By construction, $M$
contains  embedded symplectic tori with self-intersection  zero   and
the inclusion of these  tori induces the trivial  morphism on
fundamental groups.
Thus one can take  symplectic  fiber sums with arbitrarily many
(elliptically  fibered) K3 surfaces,  without changing the  fundamental
   group.  Furthermore, the fiber sums continue to be minimal by a result
of Li and Stipsicz \cite{symp:li-stip:minimal_connect_sum}.  Each such
sum increases $\chi$  by $24$ and  decreases $\sigma$ by $16$.
Therefore   $\chi+b\sigma$ can be made as small as desired
when $b>\frac32$.\end{proof}

The proof of Proposition \ref{prop:minimal_symp_man} shows that except for surface groups (of genus $g>1$), $f_{\M^{min}(G)}(2,3)=0$.   For surface groups of genus  $g>1$,
the only minimal symplectic manifolds with $f_{\M^{min}(G)}(2,3)<0$ are diffeomorphic to irrational ruled surfaces, in which case it is known that $f_{\M^{min}(G)}(2,3)=2(2-2g)$.

\medskip

Before moving on it is worthwhile to mention the consequence of  the conjectured Bogomolov-Miyaoka-Yau inequality for
symplectic manifolds to determining the  shape of $D_{\M^{min}(G)}$. Recall that the BMY conjecture states
that $\chi -3\sigma\ge 0$ for all minimal symplectic manifolds with
$K^2\geq 0$.   This gives a
lower bound for   $\chi -3\sigma$ on $\M^{min}(G)$ whenever
that $G$ is not a surface group, and hence implies that in this case $\M^{min}(G)$ contains the cone over the interval $\{1\}\times [-3,\frac32]$,
improving Proposition \ref{prop:minimal_symp_man} for non-surface groups.  It is worth  noting that  all currently known simply-connected irreducible
4--manifolds satisfy $\chi-\frac32\sigma \geq 0$.

\medskip

It is perhaps most natural to describe the domains $D_\M$ as cones on an interval contained in the unit circle and $f_\M$ as functions on these intervals.  For example $D_{\M^\infty}(G)$ corresponds to  the interval
$[-\frac\pi4,\frac\pi4]$, and $D_{\M}(G)$ corresponds to   the interval   $[\theta_G,\frac\pi4]$. However, we find it more convenient  to describe them in terms of intervals  in  $\{1\}\times \RR$ for two reasons. First,   $a\chi + b\sigma$ is not bounded  below on $\M(G)$   for $a\leq 0$.    But for $a>0$ one can divide by $a$ and minimize the 1-parameter family $\chi+b\sigma$ without losing information.
Secondly, the function of one variable
$b\mapsto f_\M(1,b)$  can easily be shown to be a piecewise linear concave function,  and   can often be explicitly described.

Thus we restrict to  the case $a>0$ (and hence to $a=1$ by normalizing) and consider the intersection of the line $\{ a=1\}$ with the domains $D_{\M^\infty}(G)$, $D_{\M}(G)$, and  $D_{\M^{min}}(G)$. Propositions \ref{prop:bound_below_M_smooth} and \ref{prop:minimal_symp_man} show that there are natural breaking points at $ b=1, $ and $b=\frac32$, corresponding to $\chi+\sigma$ and $ \chi+\frac{3}{2}\sigma$. The comments after Proposition \ref{prop:minimal_symp_man} completely compute the minimum of  $\chi+\frac{3}{2}\sigma$ on  $\M^{min}(G)$. These breaking points really do matter, as the next few calculations of the functions $f_{\M}$ over the line $a=1$ show.

\medskip

Consider first  $G=\{e\}$ the trivial group:
$$f_{\M^\infty (e)}(1,b)=\begin{cases} 2& \text {if }|b|\leq 1,\\ -\infty& \text{ otherwise.}
\end{cases}
$$
with $S^4$ the minimizer for all  $|b|\leq  1$.
By contrast,
$$f_{\M(e)}(1,b)=\begin{cases}  b+3& \text {if }|b|\leq 1,\\ -\infty& \text{if } b<-\frac{10}{3} \text{ or } b>1,\\
\text{unknown, but }\leq b+3& \text{if }-\frac{10}{3}\leq b<-1.
\end{cases}
$$
with $\CC P^2$ the minimizer for all  $|b|\leq  1$, and Stipsicz's examples \cite{stip}  treating the cases $ b<-\frac{10}{3}$.

For $\M^{min}(e)$ the domain of $ f_{M^{min}(e)}(1,b)$ includes $1\leq b\leq \frac32$. Considering   $\CP^2$,  Dolgachev  surfaces, and Stipsicz's examples yields the following:
$$f_{\M^{min}(e)}(1,b) \begin{cases} \leq b+3& \text {if } b<-1,\\
 =b+3& \text {if }|b|\leq 1,\\
=-8b +12 &\text{if }1\leq b\leq \frac{3}{2},\\
= -\infty& \text{if }b<-\frac{10}{3}\text{ or } b>\frac{3}{2}.\end{cases}
$$
Altogether, the functions yield the following graphs.

\begin{figure}[h]\label{functions}
\psfrag{b}{$b$}\psfrag{1}{$1$}\psfrag{-1}{$-1$}\psfrag{3/2}{$\frac32$}\psfrag{fM1}{$f_{\M^\infty(e)}(1,b)$}
\psfrag{fM2}{$f_{\M(e)}(1,b)$}\psfrag{fM3}{$f_{\M^{min}(e)}(1,b)$}
\includegraphics[scale=.8]{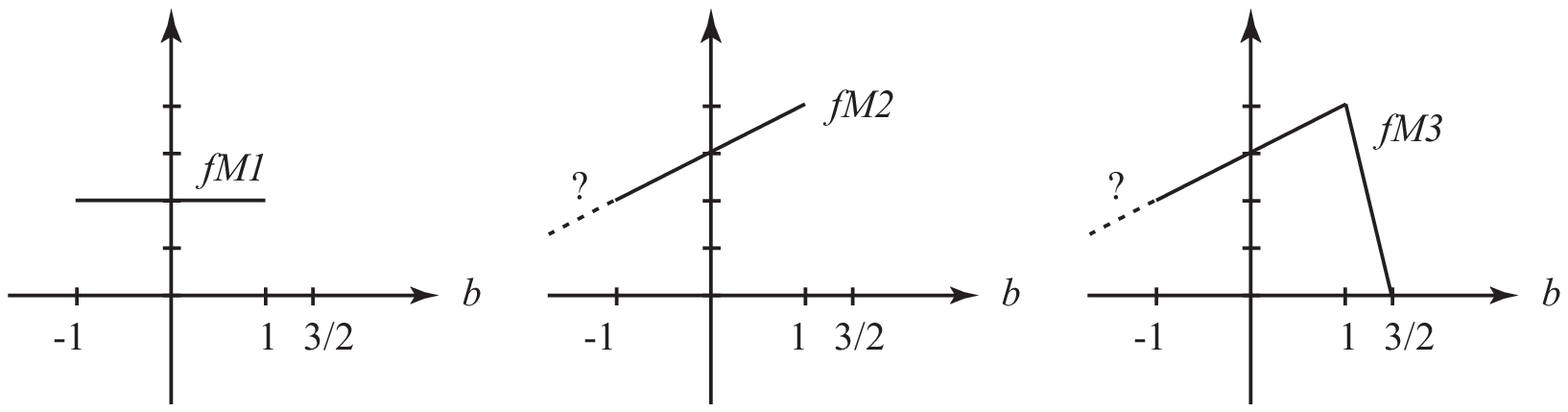}
\caption{}
\end{figure}

\medskip

Another interesting example is the  case of $G=\ZZ^6$. In \cite{group:h-w_4-man_invar_abel_grp:kirk} it is shown that  any  smooth oriented 4-manifold $M$ with fundamental group $\ZZ^6$ has $\chi(M)\ge 6$. The symplectic manifold $S_3$ described below in Section \ref{sec6} has fundamental group $\ZZ^6$,  $\chi(S_3)=6$, and $\sigma(S_3)=-2$.
Thus $f_{\M^\infty(\ZZ^6)}(1,0)$=6, and
  \begin{equation}
  \label{poss4}
f_{\M^\infty(\ZZ^6)}(1,b) \begin{cases}  \leq 6-2b& \text {if }0\leq b\leq 1,\\
=6&\text {if } b=0,\\
\leq 6+2b & \text {if }-1\leq b\leq 0, \\ -\infty& \text {if }|b|>1.
\end{cases}
\end{equation}
and hence $f_{\M^\infty(\ZZ^6)}(1,b)$ is not linear on its domain. We suspect that the inequalities in (\ref{poss4})
are equalities. This is true if and only if   $\chi(M)+\sigma(M)\ge 4$ among smooth manifolds with fundamental group $\ZZ^6$.

In (\ref{poss4}) we used the fact that  by reversing orientation shows that  $b\mapsto f_{\M^\infty(G)}(1,b)$
 is an even function.  This is not true for $\M(G)$, i.e. for symplectic manifolds, as the example  with $G=\{e\}$ above shows.

  \medskip

The domains $D_{\M^\infty(G)}$ are independent of $G$, but we do not know the answer to the question:
\begin{quest} Are the domains $D_{\M(G)}$ and $D_{\M^{\min}(G)}$ independent of $G$?
\end{quest}
The examples given above show that the functions $f_{\M(G)}$ do depend on $G$ in interesting ways.

\medskip

Motivated by the results of this section we will concentrate on minimizing $\chi$ and $\chi+\sigma$ for
the rest of the paper.

\section{Algebraic upper bounds}

We next state and prove two   theorems which give algebraically
determined bounds in terms of a presentation of $G$.

\begin{thm} \label{thm1} Let  $G$ have a   presentation with $g$
generators
$x_1,\cdots, x_g$ and $r$ relations $w_1,\cdots,w_r$. Then
   there exists a symplectic 4-manifold $M$ with $\pi_1M\cong G$,  Euler
characteristic
$\chi(M)=12(g+r+1)$, and signature $\sigma(M)=-8(g+r+1)$.
   \end{thm}

   Combining this with the  bound (\ref{hopf}) one obtains:

   \begin{cor}\label{bds3} For a finitely presented group  $G$ with  $g$
generators and $r$ relations,
   \begin{equation}\label{generalbounds}
   2-2b_1(G)+b_2(G)\leq \min_{M\in \M(G)}\chi(M)\leq  12(g+r+1).
    \end{equation}
    and
    \begin{equation}\label{generalbounds2}
   \min_{M\in \M(G)}\chi(M)+\sigma(M)\leq  4(g+r+1).
    \end{equation}
    \qed
   \end{cor}

For specific groups one can  (and we will; see  below) do better.  One
general class of groups for  which we can improve the construction of
Theorem \ref{thm1} and hence upper bound in (\ref{generalbounds}) is
treated in the following theorem. We  will show  below  that this class
includes free groups.
\begin{thm}  \label{thm2}  Let $H:F\to F$ be an orientation-preserving
diffeomorphism of an orientable surface $F$. Assume $H$ fixes a base
   point $z$. Let $G$ be the quotient of $\pi_1(F,z)$ by the
   normal subgroup generated by the words $x^{-1}H_*(x), \
x\in\pi_1(F,z)$.

     Then there exists a symplectic 4-manifold $M$ with $\pi_1M\cong G$,
Euler characteristic
$\chi(M)=12$, and signature $-8$.
   \end{thm}

\begin{proof} We prove Theorems \ref{thm1}  and \ref{thm2}
simultaneously.
The  arguments we give are derived  from Gompf's arguments
and   follow by combining them with the construction of symplectic
forms on $M\times S^1$, where $M$ is a fibered 3-manifold. The
flexibility gained by replacing Gompf's choice of $M=F\times S^1$
with a fibered manifold leads to a simplified and ultimately
smaller (as measured by the Euler characteristic) construction.

   We  begin with a discussion of how to put symplectic forms on
    4-manifolds of the form $N\times S^1$,  where $N$ is a surface bundle
     over $S^1$. This construction has its origins in Thurston's article
\cite{thurston}.

     Let $F$ be an oriented surface.  Let $H:F\to F$ be a diffeomorphism
      with at least one fixed point, and let $p:M\to S^1$ denote the
mapping torus of
      $H$, fibered over the circle with fiber $F$ and monodromy $H$.

   Let $g_0$ be a Riemannian metric on $F$, and let $g_t$ be a path of
Riemannian metrics from $g_0$ to $g_1 =H^*(g_0)$. Then $H:(F,g_0)\to
(H, g_1)$ is an isometry.

Notice that if $H$ is an isometry with respect to some metric $g_0$
then one can  take $g_t$ to be the constant path. In  this case the
volume form of $g_0$ on $F$ determines a closed 2-form $\beta$ on  $M$
whose restriction to each fiber is a volume form (i.e. a closed,
nowhere-zero, top dimensional form).

In general, we find such a 2-form as follows.  Let $\alpha_t\in
\Omega^2_F$ denote the volume form of the metric $g_t$ and  the given
orientation.
Since $H$ is an orientation-preserving diffeomorphism, the  cohomology
classes   $[\alpha_0]$  and $[\alpha_1]$ in $H^2(F;\RR)\cong \RR$ are
equal.  Hence there exists a positive smooth  function $f:[0,1]\to
(0,\infty)$ with $f(0)=1=f(1)$ so that the cohomology class
$[f(t)\alpha_t]$ is independent of $t$.  Denote   the closed,
nondegenerate 2-form $f(t)\alpha_t$ on $F$   by $\beta_t$.

Moser's stability theorem (see \cite{McDuff-Salamon}) implies that
there is a 1-parameter family  of diffeomorphisms $\psi_t:F \to F$
so that $\psi_0$ is the identity and $\psi_t^*(\beta_t)=\beta_0$. The
trace $(x,t)\mapsto (\psi_t(x),t)$
induces a  diffeomorphism $\Psi:M\to M'$, where $M'$ denotes the
mapping torus of $\psi_1\circ H$.

Let $\pi:F\times[0,1]\to F$ denote the projection to the first factor.
The  2-form $\beta$ on $F\times[0,1]$ defined by $\beta=\pi^*(\beta_0)$
   is closed. Moreover, since $(\psi_1\circ
H)^*(\beta_0)=H^*(\psi_1^*(\beta_0))=H^*(\beta_1)=\beta_0$,  $\beta$
descends to a well-defined closed 2-form on $M'$ whose restriction to
each  fiber is a  volume form. Pulling this form back to $M$ via $\Psi$
determines a closed 2-form on $M$ whose restriction  to each fiber is a
volume  form. Denote this 2-form by $\beta\in \Omega^2_M$.

Let $dt$ denote the volume form on  the base of the fibration $p:M\to
S^1$. Then $p^*(dt)$ is a 1-form on $M$. Denote by  $N $ the 4-manifold
$M\times S^1$. To distinguish it from the base of the fibration denote
the volume 1-form
on the second factor  by $ds$. Let $q_1:M\times S^1\to M$ and
$q_2:M\times S^1\to S^1$ denote  the  projections to each factor. Then
$q_2^*(ds)$ is a 1-form on $N$.

The 2-form \begin{equation}\label{defnomega}
\omega=q_1^*( \beta) + p^*(dt)\wedge q_2^*(ds)
\end{equation}
is a symplectic form on $N$.  Indeed, since  $\beta$ is closed,
$d\omega=0$, and one can check  locally that $\omega\wedge \omega$ is
nowhere zero.

If $z$ is a fixed point of   $H$, then the circle $z\times_H S^1\subset
M$ determines a torus
$T_0=(z\times_H S^1)\times S^1 \subset  M\times S^1=N$.   The
restriction of $\omega$ to this torus is a volume form; with a slight
abuse of notation it is just the form $dt\wedge ds$. Thus $T_0$ is a
symplectic torus in $N$.  Note that the self-intersection number
$T_0\cdot T_0$ in  $N$ is zero.

The fundamental  group of
$M$  is the  HNN extension of $\pi_1F$  with respect to the
automorphism induced by $H$,  i.e.
$$\pi_1 M =\langle \pi_1F, t\ | \ H_*(x)=txt^{-1} \text{ for each }
x\in \pi_1F\rangle,$$
and $\pi_1N=\pi_1M\times  \ZZ$.    Denote by $s$ the generator of the
second  factor.
Note that the Euler  characteristic and  signature of $N$ vanish.

Theorem \ref{thm2} can now  be proved, following Gompf's argument. The
group $G$ of Theorem \ref{thm2} is obtained by
taking the quotient of $\pi_1(N)$ by the normal subgroup generated by
$t$ and $s$.

Gompf's symplectic sum theorem shows that  if $E$ is a symplectic
manifold   which contains a symplectic torus $T$ with self-intersection
number zero then  the symplectic sum of $E$  and $N$,  obtained by
removing a neighborhood of the symplectic torus $T$ in $E$  and  $T_0$
in $N$ and identifying the resulting manifolds along their boundary
appropriately, then the result admits a  symplectic  structure. If,
moreover,  $\pi_1(E-T)=1$, then Van Kampen's theorem implies that the
fundamental group of the sum  $N\#_T E$ is obtained from the
fundamental  group of $N$ by killing the image of $\pi_1(T_0)$ in
$\pi_1(N)$.  Taking $E$ to be the elliptic surface $E(1)$ and $T$ a
generic fiber gives the desired symplectic  manifold $S=N\#_TE(1)$ with
$\pi_1S=G$,  $\chi(S)=12$, and $\sigma(S)=-8$.

We turn now to the  proof of Theorem \ref{thm1}.   From the
presentation  of $G$ with generators $x_1,\cdots, x_g$ and
relations $w_1,\cdots, w_r$, construct a new presentation with
$2g$ generators $x_1,y_2,\cdots, x_g,y_g$, and $g+r$ relations:
the first $g$ relations are $x_1y_2,\cdots, x_gy_g$ and the last
$r$ relations are $w_1',\cdots, w_r'$.  Here $w_i'$ is obtained
from $w_i$ by replacing every occurrence of $x_j^{-a}$ for $a>0$
with $y_j^a$ for all $j$. The relevant observation for our
purposes is that in every relation the generators appear with only
positive powers.

   Let $T=S^1\times S^1$  and define $f:S^1\times S^1\to S^1$ by
$f(e^{ia},e^{ib})=e^{i(a+b)}$. Let  $X=S^1\times\{1\}$ and
$Y=\{1\}\times S^1$. Let $D\subset T$ be a small 2-disk in the
complement of $X\cup Y$. Let $w:T\to T$ be a smooth map that  collapses
$D$ to a point and is a diffeomorphism on the complement of  $D$.
Denote by $\theta$ the 1 form  on  $T$ obtained by pulling back  the
volume form on $S^1$, $\theta=w^*(f^*(dt))$. This is a 1-form on $T$
which vanishes on  $D$,  and restricts  to a  volume 1-form on any {\em
positive monotonic} path in $T-D$, that is,  any smooth (oriented) path
   in  $T-D$ whose  composite with  $f\circ w$ wraps monotonically  (with
non-vanishing derivative)  around  $S^1$ in the positive direction.

   Let $n_i$ denote the length of the relation $w'_i$ (e.g. the length of
$x_5^3y_1y_2^2$ is $3$).  Let
$$n=1+(\sum_{i=1}^r n_i).$$

   Consider the (isometric) $\ZZ/(ng)$ action of $S^2$ generated by the
rotation $R$ about the $z$ axis by angle $2\pi/(ng)$. Let $D'$ be a
small  disc  in $S^2$ centered on the  equator (say at $(1,0,0)$) such
that its translates by
   $R$ are all pairwise disjoint.
Let $F$ be the orientable surface of genus $gn$ constructed by
removing all the translates of $D'$ by powers of $R$ and gluing in one
copy of $T-D$ along each boundary circle.   There is a corresponding
isometry $R:F\to  F$ which takes each copy of $T-D$ to the next. The
1-form $\theta$ on $T$ defines a  smooth  1-form  (which we continue to
call $\theta$) on $F$ which vanishes outside  the union of the $T-D$
and which is invariant under $R$.  Another description of this entire
construction  is to consider the $ng$-fold cyclic  branched cover of
$T$ branched over two  points in $D$ and to  pull back the 1-form
$\theta$ to the branched covering.

For convenience denote $F=A\cup B$, where $A$ is the complement of the
$ng$ discs  $R^k(D')$ in $S^2$ and $B$ is the disjoint union of the
$ng$ punctured tori.

Let $H=R^g:F\to  F$.  Thus $H$  is  an  isometry of order  $n$.   We
label the image of the curves $X$ and $Y$ in  the various copies of
$T-D$ using a double index, $X_{i,j}, Y_{i,j} ,\  i=1,\cdots,g, \
j=1,\cdots n$ labeled lexicographically. Thus $H(X_{i,j})=X_{i, j+1}$
and $ H(Y_{i,j})=Y_{i,j+1} $ (with $j$ taken modulo  $n$).  In other
words, the labeling is lifted from the $n$-fold branched cover $F\to
F/H$.

\begin{figure}
\includegraphics{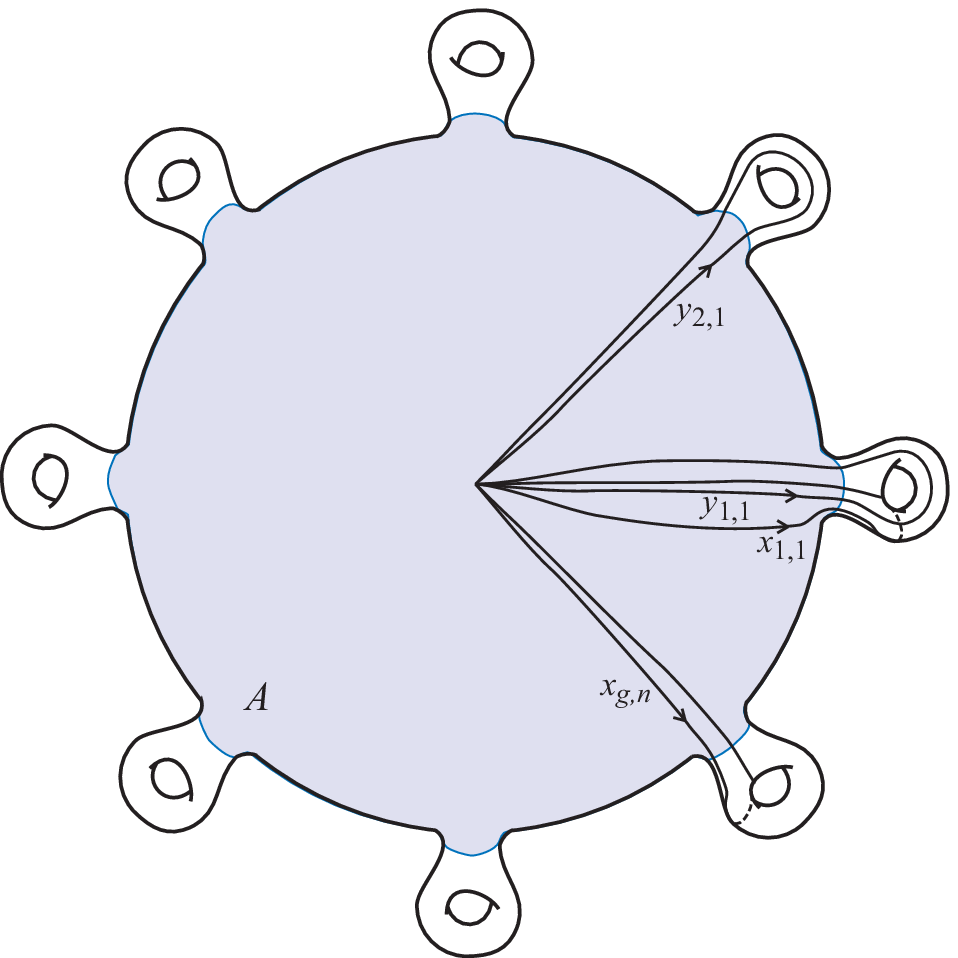}
\caption{}
\end{figure}

Join the intersection point  of $X_{i,j}$ and $Y_{i,j}$ to the north
pole $z=(0,0,1)$ along a great circle to obtain generators $x_{i,j}$
and $y_{i,j}$ of $\pi_1(F,z)$.    Thus the induced action on $\pi_1F$
is  given by $H_*(x_{i,j})=x_{i,j+1}$ and  $H_*(y_{i,j})=y_{i,j+1}$.

    To the ordered  set of    relations $w_1',  w_2', \cdots ,  w_r'$ we
assign  an  ordered set of words $\tilde{w}_1,\cdots, \tilde{w}_r$ in
the $x_{i,j}$ and $y_{i,j}$ as follows. Starting with  the first letter
    which appears in $w_1'$ replace the corresponding $x_i$ or $y_i$ by
$x_{i,1}$ or $y_{i,1}$.    For the second letter which appears in
$w_2'$ add the  second index $2$ to its  subscript, and  continue until
all the letters in $w_1'$ are replaced by doubly indexed letters in
such a  way that as  the word is read from left to right, the second
indices  increase. Then  proceed to the  second relation  $w_2'$, and
so forth. Thus  when the words $\tilde{w}_1,\cdots , \tilde{w}_r$ are
read from left to  right,  the second index in the subscripts  will
read ``$1,2,\cdots, n-1$''.

    For example, this process    converts the set  of relations
$$(w_1',w_2')=(y_2 x_1^3 x_5, y_4 y_3^2)$$
    to  $$(\tilde{w}_1,\tilde{w}_2)=(y_{2,1}x_{1,2}^3x_{5,3},
y_{4,4}y_{3,5}^2).$$

   From the $\tilde{w}_i$ one can easily construct {\em pairwise disjoint
immersed} curves $\gamma_i:S^1\to F$ for $i=1,\cdots, r$  with the
properties:
   \begin{enumerate}
\item $\gamma_i$  (connected to the north pole along a great  circle)
represents the word $\tilde{w}_i$ in $\pi_1(F,z)$.
\item The double points (if any) of $\gamma_i$ are finite,  transverse
and contained entirely in
$B$.
\item $\gamma_i$ restricts to a positive monotonic  path  in each
component (i.e.  punctured torus) of $B$.  (This is where we use the
fact that the  relations  involve  only positive powers of the  $x_i$
and  $y_i$.)

\item The curves $\gamma_i$ intersect each component (i.e. circle) of
$A\cap B$ transversely.
\end{enumerate}

The pulled back 1-form $\gamma_i^*(\theta)$ is a positive multiple
of $dt$ on  that part of $S^1$ mapped  into  the interior of  $B$
by $\gamma_i$ and is zero on  $\gamma_i\cap A$. One can find a
function $f_i$  on  $ \gamma_i^{-1}(A)$ so that $f_i$ vanishes on
the endpoint of each  arc in $ \gamma_i^{-1}(A )$ and so  that
$\gamma_i^*(\theta) + df_i$ is a volume form on  $S^1$. Since the
intersection of the union of the $\gamma_i$ with $A$ is a collection
of pairwise disjoint embedded arcs, one can  extend each $f_i$ to a
function on  $F$  which  vanishes outside a neighborhood of
$\gamma_i\cap A$ and vanishes on  $B$. Adding the values of the
$f_i$ yields  a function $f:F\to  \RR$ so  that $\gamma_i^*(\theta +
df)$ is  a  volume (i.e.  nowhere vanishing)  1-form on $S^1$ for
each $i$.

We also need $g$ extra curves, corresponding to the relations $x_iy_i,
\ i=1,\cdots, g$. Notice that $n$ was the sum of the lengths of the
relations, plus one. We use this extra bit of surface to  construct
immersions (in  fact these can be taken to be embeddings)
$\gamma_{r+k}:S^1\to F, \ k=1,\cdots, g$ corresponding to the words
$x_{1,n}y_{1,n}, \cdots , x_{g,n}y_{g,n}$.  These curves can each be
taken to  lie entirely in  one punctured torus component of $B$ and be
positive and monotonic in this component. (Alternatively, we could
have made $n$ larger and treated these relations exactly as we did with
the first type of relation. We choose this  approach since our
intention  is to  find as small a universal construction as possible.)
The 1-form $df$ vanishes on these last $g$ punctured tori  by
construction, and so $\gamma_i^*(\theta+ df)$ is a volume form for
$i=r+1,\cdots, r+g$ as well.

Since the form $\theta$ is invariant under $H$,  the pull back
$\pi_1^*(\theta)$ via the projection $ \pi_1:F\times [0,1]\to F$ is  a
closed 1-form  on $F\times [0,1]$ which determines uniquely a well
defined 1-form $\Theta$ on the mapping torus $M=F\times_H S^1$  of $H$
with the property that the restriction  of $\Theta$  to $F\times
\{0\}\subset M$ equals  $\theta$.  The function $f:F=F\times  \{0\} \to
   \RR$ extends to a  function (still called $f$) on $M$ (say by using a
cut-off  function  in the interval coordinate).
Thus we  end up with  a closed 1-form  $\Theta  +df$  on $M$ whose
restriction to  the fiber $F\times \{0\}$ pulls back   to a volume form
for each $\gamma_i:S^1\to F$.

Let $N=M\times S^1$. For small enough $\epsilon$, the form
$$\omega_\epsilon= q_1^*(\beta) + p^*(dt)\wedge q_2^*(ds) + \epsilon
q_1^*( \Theta +df )\wedge q_2^*(ds)  $$ is  a symplectic form on $N$.
For each $i=1,\cdots, r+g$ the immersed torus $T_i=\gamma_i \times S^1$
is Lagrangian with respect to $q_1^*(\beta) + p^*(dt)\wedge q_2^*(ds)$.
Since $q_1^*(\Theta+df)\wedge q_2^*(ds)$ is a volume form on
$\gamma_i\times S^1$, the  $T_i$ are symplectic with respect to
$\omega_\epsilon$ for   small  positive $\epsilon$.  The $T_i$ can  be
regularly homotoped to embeddings by a  small regular homotopy by
separating the double points of $\gamma_i$ using the parameter
transverse to the fibers in the fibration of $M$. Pushing the curve
$\gamma_i$ into a far away fiber can be used to construct a homotopy of
$T_i$ off itself. Thus the $T_i$ have self-intersection  zero.

Finally, we saw before that the ``vertical torus'' $T_0=T=z\times_H
S^1$ is symplectic with respect to  $\omega=\omega_0$; hence it
remains symplectic with respect to $\omega_\epsilon$ for  small
enough $\epsilon$.

   The fundamental  group of $N$ is generated by the $x_{i,j}$,
$y_{i,j}$, $t$, and  $s$ subject to the relations:
   $$  \prod_{i,j}[x_{i,j},y_{i,j}]=1,  t x_{i,j}t^{-1}= x_{i, j+1},
   t y_{i,j}t^{-1}= y_{i, j+1}, s\text{ is central. } $$
It follows that the quotient of $\pi_1N$ obtained by killing the
generators $s,t$, the words $\tilde{w}_i$ and $x_{i,n}y_{i,n}$ has the
presentation with  generators $x_i,\ i=1,\cdots ,g$ and relations $w_i$.

Thus to complete  the  argument we form the symplectic sum  of $N$ with
$g+r+1$ copies of the elliptic surface $E(1)$ along the symplectic tori
$T_0,T_1,\cdots,T_{r+g}$.  Summing along $T_0$ kills  $t$ and  $s$.
Summing along $T_i,\ i=1,\cdots, r$ kills $\tilde{w}_i$, and summing
along $T_{r+1},\cdots , T_{r+g}$ sets $x_{i,j}$ equal to $y_{i,j}$.
Note that this kills the commutator $[x_{i,j},y_{i,j}]$ and hence the
surface relation  disappears. A simple calculation using the
Mayer-Vietoris sequence and Novikov additivity shows that each   sum
increases $\chi$ by 12 and decreases $\sigma$ by 8, completing the
proof.
\end{proof}

Notice that the manifold $M$ constructed in the proof of Theorem
\ref{thm1}  is fibered over $S^1$ with finite order monodromy and with
two fixed points. It follows that $M$ is Seifert-fibered over a surface
$S$ of genus $g$ with  two singular fibers. If $s:M\to S$ denotes the
Seifert fibration, then  the  composite of the projection  $M\times
S^1\to M$  and  $s:M\to S$ is a singular fibration  with torus fibers.
The torus $T_0$ is one of the singular fibers. Nearby smooth fibers
form an  $n$ fold cover of $T_0$. The tori $T_i$ are products of curves
$\gamma_i$ in a section of the Seifert fibration with the  last  $S^1$
factor.

The proof of Theorem \ref{thm1} also proves the following, which is
useful  for certain classes of groups.

\begin{cor} \label{possur} Let $G$ be the quotient of a surface group
$\langle x_i,y_i\ | \ \prod_i[x_i,y_i]\rangle$ by a normal subgroup
generated by $n$ words $w_1,\cdots ,w_n$ in which the  $x_i$ and $y_i$
appear with only  positive exponents.  Then there is a closed
symplectic 4-manifold with fundamental  group $G$, Euler characteristic
$12(n+1)$ and signature $-8(n+1)$. \qed
\end{cor}

\bigskip

A very interesting question is whether the  number 12 which occurs in
Theorems
\ref{thm1} and  \ref{thm2}, and Corollaries  \ref{bds3} and
\ref{possur} can be improved.  Suppose that
$E$ is a symplectic manifold which contains a symplectic torus
$T\subset E$ such  that $T\cdot T=0$, and so  that $\pi_1(E-T)=1$.
Then if $k=\chi(E)$, the number 12 in these theorems can  be
replaced by $k$.

We can require even less: suppose that $K$ is a symplectic manifold
which contains a symplectic torus $T\subset K$ such  that $T\cdot
T=0$ and so  that $\pi_1(K-T)=\ZZ$. Let $p:T\to  K-T$ denote a  push
off of $T$ into the boundary of its tubular neighborhood.  Suppose
that the induced homomorphism  $p_*:\pi_1(T)\to \pi_1(K-T)$ is
surjective. Notice that $p_*$ contains a primitive vector in its
kernel, and so symplectically summing with $K$ can be used  just as
$E(1)$ was used in the proof.  If $\chi(K)=\ell$, then the
$12(g+r+1)$ which occurs in Theorem \ref{thm1} can  be replaced by
$\ell(g+r+2)$ or $\ell(g+r) + k$, with $k$ as in the previous
paragraph.  This is because the  first symplectic sum used in the
proof of Theorem \ref{thm1} (along $T_0$) is used to kill two
generators,  $t$ and $s$, whereas the subsequent  sums only need to
kill one generator at a time.


We  summarize these observations in the  following corollary for
completeness.

\begin{cor} \label{corvague} Let $E$ be
a closed symplectic 4-manifold which contains a symplectically
embedded torus $T$ with self-intersection zero such  that
$\pi_1(E-T)$ is trivial and with  $\chi(E)=k$. Let $K$ be a closed
symplectic 4-manifold which contains a symplectically embedded torus
$T$ with self-intersection zero such  that $\pi_1(K-T)\cong \ZZ$, $p_*:\pi_1(T)\to \pi_1(K-T)$ surjective,  and
$\chi(K)=\ell$. Then  if $G$ admits a presentation with $g$
generators and $r$ relations,
   \begin{equation}\label{generalbounds2}
    \min_{M\in \M(G)}\chi(M)\leq  k+ \ell(g+r)
    \end{equation}\qed
\end{cor}

Unfortunately, we  do not know of any ``small'' examples of $E$ or
$K$ as above.  The smallest example of such an $E$ we know is $E(1)$. The adjunction inequality (\cite{KM}) can be used to show that any such $E$ must have $\chi(E)\ge 6$. Since our constructions are based on taking fiber sums
with $E(1)$,  the smallest  example we know  of a $K$  as in Corollary \ref{corvague} has $\chi(K)=12$ (see Lemma \ref{Zexample} below).

\section{Bounds for specific classes of groups}\label{sec6}

In this section we derive better bounds for free groups, cyclic groups,
and free abelian groups  than those given in Corollary
\ref{generalbounds}.  In particular, we determine the lower bound for
certain free abelian groups and provide an example of a minimizer.

\subsection{Free groups}
   \begin{thm} \label{free}For any finitely generated free group $G$
there exists a symplectic 4-manifold $M$ with   fundamental group  $G$
and $\chi(M)=12$, $\sigma(M)=-8$.
   \end{thm}

   \begin{proof} Let $F$ be a surface of genus  $g$. Let $X_i$, $Y_i$,
$i=1,\cdots  ,g $ be a collection of embedded curves forming a standard
   symplectic basis for $H_1(F)$. Let $x_i, y_i\in \pi_1(F)$  be the
corresponding loops obtained by connecting the  $X_i,Y_i$ to  a base
point.  Take $H:F\to F$ to be the composite of Dehn twists along the
curves $Y_1,Y_2, \cdots, Y_g$. Then
   $H_*(x_i)=x_iy_i$ and $H_*(y_i)=y_i$.  It follows that the quotient of
$\pi_1(F)$ by the normal subgroup
   generated by $x^{-1}H_*(x), \ x\in \pi_1F$ is free with generators
$x_1,\cdots, x_g$.
Applying  Theorem \ref{thm2} finishes the argument.
\end{proof}

   \begin{cor}\label{corfree} Let $G$ denote the free  group on $n$
generators. Let $e=0$ if  $n$  is even  and $e=1$ if $  n$  is  odd.
Then
   \begin{equation*}
3-2n+e\leq \min_{M\in \M(G)}\chi(M)\leq  12.
    \end{equation*}
    and
     \begin{equation*}
4-2n+2e\leq \min_{M\in \M(G)}\chi(M)+\sigma(M)\leq  4.
    \end{equation*}
   \end{cor}
   \begin{proof}
   Theorem \ref{free} establishes the upper bounds.  Let $M$ be
symplectic with $\pi_1(M)\cong  G$. Notice that $\chi(M)=2-2n+ b^+(M)+
b^-(M)$ and $\chi(M)+\sigma(M)=2-2n+ 2b^+(M)$. Since $M$ is symplectic,
$b^+(M)\ge 1$. Moreover, since
   $1-b_1(M)+b^+(M)$ is even, $b^+(M)$ is even  if $n$ is  odd,  so that
for $n$ odd $b^+(M)\ge 2$.
\end{proof}

Notice that  for $G\cong \ZZ$ the upper and lower bounds in the second
formula of Corollary \ref{corfree} coincide. Thus our construction
gives a symplectic  4-manifold with fundamental group $\ZZ$ which
minimizes $\chi+\sigma$.

Kotschick \cite{kot4}  improves the  lower bound for $\min \chi$ in Corollary \ref{corfree} from $3-2n+e$ to $\frac{6}{5}(1-n)$ using the fact that $2\chi+3\sigma\ge 0$.

\subsection{Cyclic groups}

We begin with  an estimate for cyclic groups which uses  Theorem \ref{thm2}.
The argument we give is identical to  the argument given by  Gompf in  Proposition 6.4 of
\cite{symp:gompf:construct_symp_man}.

\begin{thm}[Gompf] \label{free}There exists a symplectic 4-manifold $M$ with
fundamental group $G\cong \ZZ/n$  satisfying $\chi(M)=12$ and
$\sigma(M)=-8$.
   \end{thm}

   \begin{proof} Let $F$ be a torus.    Take $H:F\to F$ to be
diffeomorphism which  induces the matrix
   $$\begin{pmatrix}  0&  1\\ -1&2-n\end{pmatrix}  $$
   on  $\ZZ^2=H_1(F)=\pi_1(F)$.
The quotient of $\pi_1(F)$ by the normal subgroup
   generated by $x^{-1}H_*(x), \ x\in \pi_1F$ is isomorphic to $\ZZ/n$,
since elementary row  and column operations transforms $H_*-I$ to the
diagonal matrix with  entries $n$  and $1$.
Applying  Theorem \ref{thm2} finishes the argument.
\end{proof}

    \begin{cor} Let $G=\ZZ/n$ for $n\ne  0$.    Then
   \begin{equation*}
3\leq \min_{M\in \M(G)}\chi(M)\leq  12.
    \end{equation*}
    and
     \begin{equation*}
   \min_{M\in \M(G)}\chi(M)+\sigma(M)=4.
    \end{equation*}
   \end{cor}
   \begin{proof}
   If $M$ is  symplectic  with $\pi_1(M)\cong\ZZ/n$, then
$\chi(M)=2+b_2(M)\ge 2+b^+(M)\ge 3$.   Moreover,
$\chi(M)+\sigma(M)=2+2b^+(M)\ge 4$. The upper bounds come  from Theorem
\ref{free}.
\end{proof}

Notice that if $M$ denotes the algebraic surface obtained from
$E(1)$ by performing two  logarithmic transformations of
multiplicity $p,q$ with $n=gcd(p,q)$, then $\pi_1(M)=\ZZ/n$,
$\chi(M)=12$, and $\sigma(M)=-8$. This shows that Theorem
\ref{free} can be improved: one can replace ``symplectic 4-manifold"
by ``\Kahler\ surface.''

   The examples of  Theorem 12 do not always minimize the Euler
characteristic. For example, there are smooth complex projective
surfaces
   with fundamental group  $\ZZ/5$ (Catanese) and $\ZZ/8$ (Reid) with
$\chi=10$.
   There are smooth complex projective surfaces
   with fundamental group  $\ZZ/2$ (Barlow and Reid) and $\ZZ/4$
(Godeaux) with $\chi=11$.
   These examples have  $\chi+\sigma=4$  \cite{BPV}.

\subsection{Free abelian groups}

We   turn  to some  calculations  and estimates  of the minimal  values
of $\chi, \chi+\sigma $ on $\M(G)$ for $G$ free abelian.

Recall first that for smooth 4-manifolds examples were constructed in
\cite{group:h-w_4-man_invar_abel_grp:kirk} which  minimize $\chi(M)$
over the class of {\em smooth} manifolds  $M$ with $\pi_1(M)=\ZZ^n$;
it was shown that the minimal Euler characteristic  for $n\ne 3,5$ is
$$2-2n+C(n,2) +\epsilon_n,$$
where   $C(n,2)$  denotes the binomial coefficient $n(n-1)/2$, and
$\epsilon_n$ is 1 if
$C(n,2)$ is odd and zero otherwise. For $n=3$ (resp. $n=5$) the minimal
Euler  characteristic is
$2$ (resp. 6).   We will  show below that for $n$  even virtually the
same result holds if we  minimize  over the class  of symplectic
4-manifolds. For $n$  odd the  situation  is less clear.

We  begin by setting  up  some notation and making some easy
observations.
Let $G=\ZZ^n$ and let $M$ be a smooth, closed 4-manifold with
$\pi_1(M)\cong G$.    Choose a map $f:M\to T^n$ inducing an
isomorphism on fundamental groups.   Since the cohomology ring $H^*(T^n)$ is an exterior
algebra on  $H^1(T^n)$, the induced map $f^*:H^2(T^n)\to H^2(M)$ is
(split) injective.  In particular $2-2n+ C(n,2)\leq \chi(M)$. Moreover,
$\chi(M)+\sigma(M)=2-2n+2b^+(M)\ge 2-2n$.

   Note that $\ZZ^n$ contains subgroups isomorphic to $\ZZ^n$ of
arbitrarily large finite index.  Since $\chi$ and $\sigma$ are
multiplicative with respect to finite covers,
\begin{equation}\label{bounds1}
   0\leq \chi(M)+\sigma(M).
\end{equation}

We turn  now to the search for  symplectic  examples  which minimize
$\chi$ and $\chi+\sigma$.

\begin{prop}\label{chi}   Any closed symplectic manifold $M$ with
$\pi_1(M)\cong\ZZ^n$  satisfies
$$  \chi(M)\ge \begin{cases} 2-2n+C(n,2)& \text{if $n\equiv 1 $  or $4$
{\rm  Mod} $8$,} \\ 3-2n+C(n,2) & \text{otherwise}\end{cases} $$
\end{prop}
and
$$ \chi(M)+\sigma(M)\equiv 0 \mod{4}.$$
\begin{proof} The cases $n=0,1,2$ are easy, so we assume that $n\ge 2$.
Suppose that $M$  is a closed symplectic  4-manifold  with
$\pi_1(M)\cong \ZZ^n$.  Then $\chi(M)=2-2n+b_2(M)$.  Since $M$ is
symplectic, $1-b_1(M)+  b^+(M)=1-n+b^+(M)$ is even. Hence
$2-2n+2b^+(M)=  \chi(M)+\sigma(M)\equiv 0 \mod{4}$.

The bound     (\ref{hopf})  (or see the paragraph  preceding Equation
\ref{bounds1})  implies that    $b_2(M)\ge C(n,2)$.  The theorem will
follow  if we can show that this bound can be improved to
$b_2(M)\ge C(n,2)+1$ when $n$ is not  congruent to $1$ or $4$ Mod $8$.

Assume that $b_2(M)=C(n,2)$.

As remarked in \cite{group:h-w_4-man_invar_abel_grp:kirk}, if
$b_2(M)=C(n,2)$, the injection  $f^*:H^2(T^n)\to H^2(M)$ is an
isomorphism. Since $H^*(T^n)$ is an exterior algebra (over $\ZZ$) on
$H^1(T^n)$, $H^2(T^n)$ has a basis for which each basis vector has cup
square zero.  This forces the intersection form of $M$ to be  even and
hence have even  rank.  Thus  $C(n,2)$ is even.

This proves that    $b_2(M)\ge C(n,2)+1$ whenever $C(n,2)$ is odd, i.e.
if $n=4k+2$ or $n=4k+3$.   (Notice that we did yet not use the  fact
that  $M$ was symplectic.)

Continue with the  assumption that $b_2(M)=C(n,2)$, so that $C(n,2)$ is
even.   Since we are assuming that $b_2(M)=C(n,2)$,  the intersection
form of $M$ is even, and  hence its signature is divisible by $8$. Thus
$\chi(M)\equiv 0\mod 4$, i.e.
\begin{equation}\label{eq2}
2-2n+C(n,2)\equiv 0 \mod{4}
\end{equation}
A simple  calculation establishes that if  $n=4k$, then  Equation
(\ref{eq2}) forces $k$ to be odd.
Similarly, if $n=4k+1$, then $k$ must be even. Thus we have shown  that
   with the possible exception of
   $n=8k+1$ and $n=8k$, a symplectic 4-manifold  $M$ with  $\pi_1(M)\cong
\ZZ^n$ must have
   $b_2(M)\ge C(n,2)+1$, and so $\chi(M)\ge 2-2n+C(n,2)+1$.

\end{proof}

In \cite{group:h-w_4-man_invar_abel_grp:kirk} it was shown that there
exist smooth  closed 4-manifolds $X_n$ with $\pi_1(X_n)\cong\ZZ^n$ and
$\chi(X_n)=2-2n+C(n,2)$ for any $n>5$ with $C(n,2)$ even.
It follows from Theorem \ref{chi} that $X_n$ cannot admit a symplectic
structure when
$n=8k $ or $n=8k+5$.  For these examples, $b_1(X_n)$ is even,
$\sigma(X_n)=0$ and $2\chi(X_n)+3\sigma(X_n)\ge0$.

As explained in \cite{group:h-w_4-man_invar_abel_grp:kirk}, the cases $\ZZ^3$ and $\ZZ^5$ are exceptional. The intersection form of any smooth manifold $M$ with fundamental group $\ZZ^3$ has a 3-dimensional metabolizer, hence $b^+(M)\ge 3$ and $b^-(M)\ge 3$. If $M$ is symplectic then $b^+(M)$ is even, hence at least $4$. Thus $\chi(M)\ge 3$ and $\chi(M)+\sigma(M)\ge 4$. Similarly,  the intersection form of any smooth manifold $M$ with fundamental group $\ZZ^5$ has a 7-dimensional metabolizer. If $M$ is symplectic this implies $\chi(M)\ge7$ and $\chi(M)+\sigma(M)\ge 8$.

\bigskip

We next look  at upper bounds. As a first estimate, since $\ZZ^n$ has a
presentation with  $n$ generators and $C(n,2)$ relations, Theorem
\ref{thm1} and  Proposition \ref{chi} give the estimates
$$\tfrac{1}{2} (n^2-5n+4)\leq  \min_{M\in \M(\ZZ^n)}\chi(M)\leq
6(n^2+n+2).$$
Thus we see that $ \min_{M\in \M(\ZZ^n)}\chi(M)$ grows quadratically in
$n$, with leading coefficient between $\tfrac{1}{2}$ and $6$. It
follows from the calculations below that
restricting to $n$ even, $ \min_{M\in \M(\ZZ^{2n})}\chi(M)\sim
\tfrac{1}{2}(2n)^2$, and we will give evidence that the restriction to
even rank is unnecessary.

\bigskip
For each  integer $g\ge 0$ let $F_g$ denote the  surface of  genus $g$.
   Let $S_g=\text{Sym}^2(F_g)$.

\begin{prop}\label{Sg} The space $S_g$ is  a compact \Kahler\ manifold,
and in particular is  symplectic. Moreover,
$\pi_1(S_g)=\ZZ^{2g}$,  $H^2(S_g)= \ZZ^{C(2g,2)+1}$ so that
$\chi(S_g)=3-2(2g)+C(2g,2)$, and  $\sigma(S_g)=1-g$.
\end{prop}
\begin{proof}[Sketch of proof]
The fact that $S_g$ admits a complex structure comes from the fact that
the $\ZZ/2$ action on $F_g\times F_g$ defines a branched cover of
$S_g$. The fundamental group is computed using Van-Kampen's theorem
splitting $S_g$ along the circle bundle over the branch set; note that
the two sets of generators of $\pi_1(F_g\times F_g)\cong
\pi_1(F_g)\times \pi_1(F_g)$ commute, and are identified in
$\pi_1(S_g)$. Since $b_1(S_g)$ is even, $S_g$ is \Kahler\
\cite{sw:4man_kirby_calc}.

The Riemann-Hurwitz formula computes $\chi(S_g)$ and   with the
universal coefficient theorem this implies the computation for
$H^2(S_g)$. Computing the signature is a bit more involved; the most
straightforward way to do this is to use the transfer (with $\RR$ coefficients) to observe that
the induced map $H^2(S_g)\to H^2(F_g\times F_g)$   is injective with
image the $\ZZ/2$-invariant classes, and to compute the intersection
form directly by restricting the intersection form of $F_g\times F_g$.

We refer to \cite{sym:mcdonald:symproduct_of_surface} for details.
\end{proof}

The following corollary computes the minimal Euler  characteristic for
most $\ZZ^{2g}$.

\begin{cor}\label{evenZn} Let $G=\ZZ^{2g}$.   Then \hfill
\begin{enumerate}
\item If $g\equiv 0,1$ or $3$ {\rm  Mod} $4$, then
   \begin{equation*}
    \min_{M\in \M(G)}\chi(M)=
3-4g+C(2g,2),
   \end{equation*}
   with minimizer $S_g$.
   \item  If $g\equiv  2$ {\rm  Mod}  $4$, then
    \begin{equation*}
0\leq  \min_{M\in \M(G)}\chi(M)-\big(2-4g+C(2g,2)\big)\leq 1,
   \end{equation*}

   \item $\displaystyle
   0\leq\min_{M\in \M(G)}\chi(M)+\sigma(M)\leq 4-5g+C(2g,2).
$
\end{enumerate}
\end{cor}
\begin{proof}  The examples $S_g$ of Proposition \ref{Sg}  provide the
upper bounds. Proposition
\ref{chi}  shows that when $g\equiv 0,1$ or $3$ {\rm  Mod} $4$, the
$S_g$ give the smallest possible $\chi$. When $g\equiv 2$ {\rm  Mod}
$4$, the lower  bound of Proposition \ref{chi} differs by  one from
$\chi(S_g)$.   The third assertion comes from Equation \ref{bounds1}
and Proposition \ref{Sg}.

\end{proof}

Corollary \ref{evenZn} does not  answer the question of whether $S_g$
minimizes $\chi$ on $\M(\ZZ^{2g}) $ when $g=4m+2$. In  fact it does not
for $m=0$:
    $S_2=T^4\#\overline{\CC P}^2$.  But the 4-torus $T^4$ is symplectic
and
$$0=\chi(T^4)<\chi(S_2)=1.$$
We do not know whether $S_{4m+2}$ minimizes  $\chi:\M(\ZZ^{8m+4})\to
\ZZ$ for $m>0$.

Note that $S_2$ does minimize $\chi+\sigma$. In fact, $S_0,S_1$ and
$S_2$    minimize
$\chi+\sigma$ for symplectic manifolds and  $G= 0, \ZZ^2, \ZZ^4$. The
first unknown  case is
$S_3$, with  $\chi(S_3)+\sigma(S_3)=4$.  Hence either $S_3$ minimizes
$\chi  +\sigma$ among  symplectic 4-manifolds with  fundamental  group
$\ZZ^6$ or else (since $\chi+\sigma\equiv 0$ Mod 4)
there is a symplectic 4-manifold $X$  with $\pi_1(X)=\ZZ^6$ and
$b^+(X)=5$.

\bigskip
   Free abelian groups of odd rank pose a greater challenge. For $G=
\ZZ$, we know that any symplectic 4-manifold $M$ with $\pi_1(M)\cong
\ZZ$ has $b^+(M)$ even and greater than zero, thus $\chi(M)=b_2(M)\ge
2$.  On the other hand, Theorem \ref{free} constructs a symplectic
4-manifold with $\pi_1(M)\cong \ZZ$ with $\chi(M)=12$ and
$\sigma(M)=-8$. At the moment this is the smallest example known to the
authors of a symplectic 4-manifold with fundamental group $\ZZ$ (see
Theorem \ref{free}).  Thus
   \begin{equation}\label{zbound}
    2\leq  \min_{M\in \M(\ZZ)}\chi(M)\leq 12.
   \end{equation}
   Note that the lower bound was derived using only the fact that $M$ is
an almost complex manifold rather than the stronger assumption that $M$
is symplectic.

   This example does minimize $\chi+\sigma$. Indeed, since $b^+(M)$ is
even and greater than zero for a symplectic 4-manifold with fundamental
group $\ZZ$, it follows that $\chi(M)+\sigma(M)=2b^+(M)\ge 4$. The
example of Theorem \ref{free} with $\pi_1(M)\cong\ZZ$ has
$\chi(M)+\sigma(M)=4$, so
   \begin{equation}
    \min_{M\in \M(\ZZ)}\chi(M)+\sigma(M)=4.
   \end{equation}

   \medskip
We turn to  the case $G\cong \ZZ^3$.

Consider the four-torus $X=T^2\times T^2$ with the product
symplectic structure. Its fundamental group is $\ZZ^4$ generated by
the coordinate circles; call these generators $a,b,c,d$. The Euler
characteristic of $X$ is $0$ and $\sigma(X)=0$. The symplectic torus
$T_0=p\times T^2$  has fundamental group generated  by $c$ and $d$
and self-intersection $0$.  We can use the manifold in the next
lemma to kill one of the generators $c$ or $d$.

\begin{lem}\label{Zexample} There exists a symplectic 4-manifold $K$
with $\pi_1(K)\cong \ZZ$ which contains a symplectically embedded torus
$T$ with self-intersection zero such that
\begin{enumerate}
\item $\chi(K)=12$ and $\sigma(K)=-8$.
\item $\pi_1(K-T)\cong \ZZ$ and the map induced by inclusion
$\pi_1(K-T)\to \pi_1(K)$ is an isomorphism.
\end{enumerate}
   \end{lem}
\begin{proof} Let $K$ be the symplectic 4-manifold with $\pi_1(K)\cong
\ZZ$ constructed in Theorem \ref{free}.  The construction of $K$ was
the following. First a fibered 3-manifold $M$ with fiber a torus $F$ is
constructed as the mapping torus of the Dehn twist $H:F\to F$ on the
torus along the second curve $y$ of a symplectic basis $\{x,y\}$ of
$\pi_1(F)$.   Thus $\pi_1(M)=\langle x,y,t\ | \ [x,y], txt^{-1}=xy,
tyt^{-1}=y\rangle$ and letting $N=M\times S^1$,
$\pi_1(N)= \langle x,y,t, s \ |\  [x,y], txt^{-1}=xy, tyt^{-1}=y, s
\text{ central }\rangle$.
Then $N$ contains a symplectic form $\omega$ (see Equation
(\ref{defnomega})) for which the torus $T_0=t\times s$ is symplectic,
and  taking the symplectic  fiber sum of  $N$ with $E(1)$ along $T_0$
yields $K$. Since $\pi_1(N-T_0)\to  \pi_1(N)$ is surjective  and
$\pi_1(E(1)-T_0')=1$, where $T_0'$ is the elliptic fiber in $E(1)$
along  which the symplectic sum is  taken, it follows that  $\pi_1(K)$
is infinite cyclic, generated by $x$.

   Let $T$ denote the embedded torus in $N$ given by $ x\times s$. More
precisely, choose an embedded curve $\gamma$ freely homotopic to $x$ in
the fiber $F$ which avoids the base point. Then $T=\gamma\times
S^1\subset M\times S^1$ is a torus and the morphism induced by
inclusion takes the two generators of $\pi_1(T)$ to $x$ and $s$.  From
Equation (\ref{defnomega}) one sees that $T$ is Lagrangian in $N$.
Notice also  that $T$ is disjoint from $T_0=t\times s$ since $\gamma$
avoids the base point of $F$.  Also notice that $T$ has self
intersection zero since $\gamma$ can pushed off itself.

Since $x$ is non-zero in $H_1(M)$, $T$ is non-zero in $H_2(N)$ by
the K\"unneth theorem.  It follows by a standard argument (see e.g.
Lemma 1.6 of \cite{symp:gompf:construct_symp_man}) That $\omega$ can
be perturbed by an arbitrarily small amount so that $T$ is
symplectic with respect to the   resulting symplectic form
$\omega'$. If the perturbation is taken very small, $T_0$ remains
symplectic. Gompf shows furthermore that a symplectic structure on
the fiber sum $K=N\#_{T_0=f}E(1)$ can be chosen so that $T$ remains
symplectic in $K$.

Thus $T\subset K$ is a symplectic torus with self-intersection zero
for which the induced map on fundamental groups is the map $\ZZ
x\oplus \ZZ s \to \ZZ x$, i.e. $x\mapsto x$, $s\mapsto 1$.  To
compute $\pi_1(K-T)$, first notice that $N-T=(M-\gamma)\times S^1$.
Since $\gamma$ is a curve in the fiber of the fibration $M\to S^1$
(representing $x$), it follows that   $M-\gamma$ is obtained from
$F\times[0,1]$ by gluing the ends along an  annulus, namely the
annulus in the torus $F$ complementary to $\gamma$.  Thus
$\pi_1(M-\gamma)=\langle x,y,t\ | \ [x,y], txt^{-1}=xy\rangle$. It
follows that $\pi_1(N-T)=\langle x,y,t,s\ | \ [x,y], txt^{-1}=xy, s
\text{ central }\rangle$. Since $K-T$ is the fiber sum of $N-T$ with
$E(1)$ along $T_0$,
$$\pi_1(K-T)=\langle x,y,t,s\ | \ [x,y], txt^{-1}=xy, s \text{ central
}, s=1,t=1\rangle=\ZZ x.$$

\end{proof}

Fiber sum $K$ to $T^4$ along using $T$ in $K$ and $T_0$ in $T^4$,
$$L=T^4\#_{T=T_0}K,$$
identifying $x$ with $c$ and $s$ with $d$.  This, in effect, kills
$d$ without introducing any new relations, giving a symplectic
manifold with fundamental group $\BZ^3$, with $\chi(L)=12$ and $\sigma(L)=-8$. Together with the remarks after Proposition \ref{chi} this implies the following.

\begin{prop}
There exists a symplectic $4$-manifold $L$ with fundamental group
$G=\BZ^3$ satisfying $\chi(L)=12$ and $\sigma(L)=-8$.    Hence
    $$ 3\leq \min_{M\in \cM(\ZZ^3)} \chi(M)\leq 12.$$
    Moreover,
    $$  \min_{M\in \cM(\ZZ^3)} \chi(M)+\sigma(M)=4.$$\qed
\end{prop}

Finally, we treat the  case of odd rank free abelian groups.

\begin{thm} \label{thmodd} There exists a symplectic 4-manifold $M$ with $\pi_1(M)\cong\ZZ^{2n-1}$ such that
$\chi(M)=15-5n+ 2n^2$ and  $\sigma(M)=-7-n.$
\end{thm}

Theorem  \ref{thmodd} gives the bound
$$\min_{M\in\M(\ZZ^{2n-1})}\chi(M) -\big( 2-2(2n-1)+C(2n-1,2)\big) \leq
2n+10.$$ In other words,  the difference between the lower bound of Equation
(\ref{hopf}) and the examples constructed here grows linearly  with the
rank. This is in contrast with the the examples of even rank free abelian
groups: that difference is always a constant.   On the other  hand, it is an
improvement over the general construction of Theorem  \ref{thm1}, whose
difference grows quadratically in $n$.

The proof of Theorem \ref{thmodd} depends on finding a suitable symplectic  form on  the (Kahler) manifold  $S_g=$Sym$^2(F_g)$ for which we can identify  certain tori as Lagrangian.
The main  technical result needed is the  following proposition,  whose proof was suggested to us by R. Gompf.

\begin{prop} \label{pushdownbaby} Let  $\pi:F_g\times  F_g\to  S_g$ denote  the regular 2-fold branced cover  corresponding  to  the $\ZZ/2$  action $(x,y)\mapsto (y,x)$  on $F_g\times F_g$ with fixed submanifold $B=\{(x,x)\ | \ x\in F_g\}$.  Let $\omega_F$ be a fixed symplectic form on  the surface $F$ and let $\omega=\omega_F\oplus \omega_F\in \Omega^2(F_g\times F_g)$ be the
 $\ZZ/2$-equivariant symplectic form on the product.

   Then there exists a symplectic form $\omega'\in \Omega^2(S_g)$ so that the pullback $\pi^*(\omega')$ agrees with $\omega$ outside a small tubular  neighborhood  of $B$.
\end{prop}
\begin{proof}

We first show that that in any neighborhood of $B$ in $F\times F$  one can  find a tubular neighborhood $N$ of $B$  which admits a semi-free Hamiltonian $S^1$ action with  fixed set $B$, such that the $\ZZ/2$ action $(x,y)\mapsto(y,x)$ embeds in the  $S^1$ action  as multiplication by $-1$. The Hamiltonian function $\mu:N\to  [0,\epsilon)$ satisfies $\mu^{-1}(0)=B$.
This is a standard fact in symplectic topology; we include a proof for the benefit of the reader.

Fix a Riemannian metric on $F\times F$ and let $P\to B$ be the principal $SO(2)=U(1)$ bundle associated to the normal bundle of $B$, i.e. $P$ is the bundle with $c_1(P)= 2-2g$. Then  $P$ admits a free $\ZZ/2$ action commuting with   the $U(1)$ action,  namely multiplication by  $-1\in U(1)$. Let $E=P\times_{U(1)}D^2\to B$ be the associated disc bundle.   Note that $E$ is diffeomorphic to the normal disc bundle $\nu$ of $B\subset F\times F$. Moreover, one can  choose the  diffeomorphism $E\cong \nu$ equivariant with respect to the $\ZZ/2$ action on $E$ and the linearization of the $\ZZ/2$ action $(x,y)\mapsto(y,x)$ near $B$ in  $F\times F$.

The symplectic form  $\omega$ on $F\times F$ restricts to a symplectic form $\omega_B$ on $B$.   This form  extends to an $S^1$-equivariant  symplectic form on $E$ with corresponding Hamiltonian function $\mu:E\to [0,1]$, so that $\mu^{-1}(0)=B$ (see \cite[page 155]{McDuff-Salamon}).

Since $E$ and $\nu$ are equivariantly diffeomorphic symplectic bundles and restrict to the  same symplectic form  $\omega_B$ on  $B$, Weinstein's symplectic tubular neighborhood theorem
(see \cite{Weinstein} and  \cite[page 98]{McDuff-Salamon}) implies that there is a $\ZZ/2$-equivariant  symplectomorphism from a neighborhood of  the zero section in $E$ and a neighborhood of $B$ in $F\times F$.  Since any neighborhood of the zero section in $E$ contains a smaller neighborhood of the form $E_\epsilon=\mu^{-1}([0,\epsilon))$, pulling  back $\mu$ and the $U(1)$ action  via the symplectomorphism  restricted  to $E_\epsilon$ gives the  desired neighborhood $N$,  Hamiltonian $\mu$, and   corresponding Hamiltonian $S^1$ action.

 Denote the quotient of $N-B$ by  the $\ZZ/2$ action by  $U$. Thus $U$  is endowed with  the quotient symplectic stucture (since $\ZZ/2$ acts freely and symplectically on $N-B$ with quotient $U$) and admits a free Hamiltonian $S^1(=S^1/(\ZZ/2))$ action  with Hamiltonian  $\bar{\mu}:U\to  (0,\epsilon)$.

 Symplectic cutting $U$  at $\epsilon/2$ (see \cite{lerman}) yields a  symplectic manifold $\bar{N}$  diffeomorphic to
 the tubular neighborhood of the  branch set $\bar{B}\subset S_g$.  The symplectic structure on $U$  is the restriction to $U$ of the   symplectic structure on $S_g-\bar{B}$    (pushed down from the equivariant symplectic  structure on  $F\times F-B$.)  Since symplectic cutting preserves the symplectic structure away from the  cut locus it follows that $S_g$ admits a symplectic form $\omega'$ whose restriction
 to  $S_g-\bar{N}$ pulls back to the restriction of $\omega$ to $F\times F-N$.
\end{proof}

Notice that the proof of Proposition \ref{pushdownbaby}  applies equally well to any regular branched cover $X\to  Y=X/G$ with connected,  symplectic  branch   manifold $B\subset X$ and $G$-equivariant symplectic form $\omega$ on $X$.

\medskip

\begin{proof}[Proof of Theorem \ref{thmodd}]  Let $F$ be a closed surface  of genus $g$ with a symplectic form  $\omega_F$. Let $\gamma_1$  and $\gamma_2$  be disjointly   embedded curves in $F$ representing different vectors in a symplectic basis for $H_1(F)$.   Then $T=\gamma_1\times \gamma_2$ is a Lagrangian  torus in  $F\times F$.    Since $\gamma_1$ and  $\gamma_2$ are disjoint the composite $T\subset F\times F\to S_g$ is also  an embedding.  Proposition \ref{pushdownbaby}  implies that this torus (which we continue to denote  $T$) in $S_g$ is Lagrangian  with  respect  to a suitable symplectic  form on $S_g$.       The torus $T\subset S_g$ represents a  non-trivial  homology class in $H_2(S_g)$ since  its transfer $\tau([F])\in H_2(F\times F)$ is nonzero (it equals $\gamma_1\times\gamma_2+ \gamma_1\times\gamma_2$)  by the Kunneth theorem.

Thus the  symplectic form on $S_g$ can be perturbed slightly so that $T\subset S_g$ is symplectic. Taking the symplectic fiber sum of $S_g$ with the manifold $K$ constructed in Lemma \ref{Zexample}
so that  yields a symplectic manifold $M$ whose fundamental group is the quotient of $\pi_1(S_g)=\ZZ^{2g}$ by the subgroup generated by $[\gamma_1]$, i.e. $\pi_1(M)=\ZZ^{2g-1}$, and such that $\chi(M)=\chi(S_g)+12$, $\sigma(M)=\sigma(S_g)-8$.   The calculations of Proposition \ref{Sg} finish the proof.
\end{proof}

\subsection{Other abelian groups}
In Section 6 of \cite{symp:gompf:construct_symp_man} (Propositions 6.4 and 6.6), Gompf
explores the geography of symplectic 4-manifolds with certain abelian fundamental groups
constructed by symplectically summing torus bundles with $E(1)$. For completeness
we state   his results in our terminology.
\begin{thm}[Gompf] \label{gpf}\hfill
\begin{enumerate}
\item If $G$ is the direct sum of up to three cyclic groups, except $\ZZ\oplus\ZZ\oplus \ZZ$, or if
$G=\ZZ\oplus\ZZ\oplus\ZZ/k\oplus \ZZ/\ell$ with $k,\ell\ne 0$, then there is a symplectic 4-manifold $M$ with $\pi_1(M)=G$, $\chi(M)=12$ and $\chi(M)+\sigma(M)=4$.

\item  If $G$ is   $\ZZ\oplus\ZZ/k\oplus \ZZ/\ell\oplus \ZZ/n$, or if
$G=\ZZ\oplus\ZZ\oplus\ZZ\oplus \ZZ/k$ with $k,\ell, n\ne 0$,  then there is a symplectic 4-manifold $M$
with $\pi_1(M)=G$, $\chi(M)=24$ and $\chi(M)+\sigma(M)=8$.
\end{enumerate}
\end{thm}
Note that these computations include the computations we gave for cyclic groups in the previous subsections.  Using the same arguments as in the previous subsections, the first statement in Theorem \ref{gpf} has the following consequences:
\begin{enumerate}
\item If $G=\ZZ/k\oplus \ZZ/\ell\oplus \ZZ/n$ with $k,\ell, n\ne 0$, then
$$3\leq\inf_{\M(G)}\chi(M)\leq 12 \text{ and } \inf_{\M(G)}\chi(M)+\sigma(M)=4.$$
\item If $G=\ZZ/k\oplus \ZZ/\ell\oplus \ZZ$ with $k,\ell\ne 0$, then
$$2\leq\inf_{\M(G)}\chi(M)\leq 12 \text{ and } \inf_{\M(G)}\chi(M)+\sigma(M)=4.$$
\item If $G=\ZZ/k\oplus \ZZ^2$ with $k\ne 0$, then
$$0\leq\inf_{\M(G)}\chi(M)\leq 12 \text{ and } \inf_{\M(G)}\chi(M)+\sigma(M)=0 \text{ or }4.$$
\end{enumerate}
Corresponding (but weaker) bounds can be derived from the second statement of Theorem \ref{gpf}.

Gompf also gives examples of relatively small symplectic 4-manifolds with other (non-abelian) fundamental groups. We refer the interested reader to his beautiful article \cite{symp:gompf:construct_symp_man}.

\section{Some Final Remarks}

We end with a small discussion about some difficult issues
surrounding minimizers of $\chi$. The 4-dimensional Poincar\'e
conjecture can be rephrased by saying that any  simply connected
topological (resp. smooth) 4-manifold  with minimal Euler
characteristic  is homeomorphic (resp. diffeomorphic) to the
4-sphere. In other words, if one minimizes the Euler characteristic
$\chi$ on the  class of  simply connected 4-manifolds, the minimizer
is unique.   Freedman's theorem \cite{freedman} proves the
Poincar\'e conjecture for topological  manifolds, and the smooth
question is  one of the outstanding problems in 4-dimensional topology.

New wrinkles appear in the symplectic  case.  For example $\CP^2$
minimizes the Euler characteristic among simply-connected symplectic
4-manifolds,  and
Freedman's theorem implies any two minimizers are homeomorphic. One
might call the problem of whether any two simply connected
symplectic 4-manifolds with $\chi=3$ are diffeomorphic (or
symplectomorphic) the ``symplectic Poincar\'e conjecture''.  A
counterexample would involve finding a simply-connected symplectic
$4$-manifold $(M,\omega)$ having $\chi(M)=3$ and $K_M\cdot
[\omega]>0$ (c.f. \cite{LL} or \cite{symp:app_gen_wall_cross}).  The
question of whether a simply connected symplectic manifold with
$\chi=3$ is diffeomorphic or symplectomorphic to $\CP^2$ is
unresolved, but there has been much recent progress in the direction
of a counterexample.  Starting with
\cite{sw:park:sc_symp_4_man_chi_is_10} and expanded upon in
\cite{symp:ozsvath:On_Park_manifolds,sw:stip-n-szabo:exotic_smooth_structure_chi_is_9,
sym:Ron-n-Ron:Double_Node,sw:park:exotic_smooth_chi_is_7}, new
examples were constructed of irreducible smooth $4$-manifolds
homeomorphic but not
diffeomorphic to $\CP^2\#n\overline{\CP^2}$ for $n=5,6,7,8$.  However,
for $n=5$ the examples are not symplectic.

All attempts to change the diffeomorphism type of known minimizers
without changing their fundamental group  seem to fail,  suggesting
that minimizers of $\chi:\M(G)\ra \BZ$ are somehow special.
But to conjecture that a symplectic minimizer of $\chi:\M(G)\ra \BZ$  is
unique up to diffeomorphism, however, is simply incorrect.  For
example, $S^2\times T^2$ and the nontrivial
$S^2$-bundle over $T^2$ both have fundamental group $\BZ^2$.  Yet the
search for other examples with $G=\BZ^2$ seems futile.  It is certainly
easy to build homology $T^2\times S^2$ symplectic manifolds:  let $Y$
be zero surgery on a fibered knot in $S^3$ and take $Y\times S^1$.  The
only example from this extensive list that has fundamental group
$\BZ^2$ is when the knot is the unknot, i.e., when $Y\times S^1$ is
diffeomorphic to $T^2\times S^2$.  The key difference between
$S^2\times T^2$ and the nontrivial $S^2$-bundle over $T^2$ is that the
first is spin and the second is not.  So minimizers of $\chi:\M(G)\ra
\BZ$ can have different intersection forms of the same rank.  This leads
us to make, possibly out of ignorance, the following conjecture:

\begin{conj}  Let   $M$ be a symplectic  4-manifold with $\pi_1(M)\cong
   G$ which  minimizes   $\chi:\M(G)\ra \ZZ$.  Let $Q_M$ denote the
intersection form of $M$.  Then any  other symplectic manifold   with
intersection form $Q_M$   which also  minimizes
$\chi:\M(G)\ra \ZZ$ is diffeomorphic to  $M$. \end{conj}

We offer this conjecture merely as a new twist on an old theme in
4--manifold theory, namely, describing conditions under which
4--manifolds are possibly unique.  A weaker  conjecture would be to
let $Q_M$ denote the equivariant (i.e. $\ZZ[G]$)  intersection form of
$M$. A counterexample to this conjecture would also be interesting.
A good place to start is to find another minimizer of $\M(\BZ^6)$ which
is not diffeomorphic to $S_3$.   Notice that
any  minimizer of $\chi$ is necessarily minimal. If $G$ is not  a free
product then  any  minimizer of $\chi$ is
irreducible.

\medskip

Suppose instead that one looks for  minima of $\chi + \sigma$ on
$\M(e)$. Then minimizers
are not unique: for example  $\CP^2\# n\overline{\CP}^2$ are
minimizers in $\M(e)$. These examples indicate that to go beyond
excessively general observations one may have to restrict further
the class of manifolds, e.g. irreducible manifolds. Even then
minimizers are not unique (up to diffeomorphism, for example).
Indeed there are examples mentioned above of irreducible, symplectic
4-manifolds
homeomorphic but not diffeomorphic to $\CP^2\#n \overline{\CP}^2$
for $n=6$ (c.f.
\cite{sw:stip-n-szabo:exotic_smooth_structure_chi_is_9}).

\medskip

We end this article with   remarks about improving our bounds.

What is missing in our results is a method for increasing the lower
bounds of
   $\min_{M\in \M(G)}\chi(M)$ which uses the fact that $M$ is  symplectic
in a non-trivial  way.  The lower bounds given in  the  present
article  are  obtained by  combining the lower bounds valid  for all
4-dimensional Poincar\'e complexes (e.g.  Equation (\ref{hopf}))
   with two  simple facts which hold for symplectic manifolds: $b^+(M)\ge
1$ and  $1-b_1(M)+b^+(M)$ is even. This second fact   depends only the
existence of an almost complex structure.  Our calculations show that  for $G=\ZZ^{2g}$, the difference
$$\min_{\M(\ZZ^{2g})}\chi(M)-\min_{\M^\infty(\ZZ^{2g})}\chi(M) $$
equals zero or one.  On the other hand, a recent article of Kotschick \cite{kot4} shows that for $G_k$ the free group on $k$ generators, the difference
$$\min_{\M(G_k)}\chi(M)-\min_{\M^\infty(G_k)}\chi(M) $$  gets arbitrarily large as $k$ goes to infinity.  Thus any improvement of the  lower
bounds which uses
the symplectic  structure in a deeper way will have to take  these
kinds of examples into  account.

As explained at the end of Section 4, improving our upper bounds
requires that we find a symplectic 4-manifold $K$ with     $\chi(K)<12$
which contains a symplectically embedded torus $T$ of self-intersection number zero with
$\pi_1(K-T)\cong \ZZ$ or $\pi_1(K-T)=1$. We have not found any such
manifold, and might conjecture that one does not exist. It  is  not hard to show  that any
such  $K$ must satisfy $\chi(K)\ge 6$.







\end{document}